\documentclass[12pt]{amsart}
\usepackage{amsmath}
\usepackage{amssymb}
\usepackage{latexsym}

\newcommand{\C}{\mathbb {C}}     
\newcommand{\halmos}{\rule{5pt}{5pt}}

\numberwithin{equation}{section}

\newtheorem{prop}{\bf Proposition}[section]
\newtheorem{thm}[prop]{\bf Theorem}

\newtheorem{conj}{\bf Conjecture}

\newenvironment{prf}{\noindent{\em {\normalsize P}\normalsize{roof.}}  
\normalsize\hskip 5pt}{\hfill\halmos}

\begin{document}
\title[Variants of $q$-hypergeometric equation]
{Variants of $q$-hypergeometric equation}
\author[N.~Hatano]{Naoya Hatano}
\address{Address of N.H.,~R.M.,~T.S., Department of Mathematics, Faculty of Science and Engineering, Chuo University, 1-13-27 Kasuga, Bunkyo-ku, Tokyo 112-8551, Japan}
\author[R.~Matsunawa]{Ryuya Matsunawa}
\author[T.~Sato]{Tomoki Sato}
\author[K.~Takemura]{Kouichi Takemura}
\address{Address of K.T., Department of Mathematics, Ochanomizu University, 2-1-1 Otsuka, Bunkyo-ku, Tokyo 112-8610, Japan}
\email{takemura.kouichi@ocha.ac.jp}
\begin{abstract}
We introduce two variants of $q$-hypergeometric equation.
We obtain several explicit solutions of variants of $q$-hypergeometric equation.
We show that a variant of $q$-hypergeometric equation can be obtained by a restriction of $q$-Appell equation of two variables. 
\end{abstract}
\maketitle

\section{Introduction}

One of the most important special functions is the Gauss hypergeometric function defined by
\begin{align*}
& \: _2 F _1 (\alpha ,\beta ;\gamma ;z)=1+ \frac{\alpha \beta  }{ \gamma } z + \frac{\alpha (\alpha +1) \beta (\beta +1 ) }{2! \: \gamma (\gamma +1)} z^2 + \dots +\frac{(\alpha )_n (\beta )_n }{n!(\gamma )_n } z^n + \dots ,
\end{align*}
where $(\lambda )_n= \lambda (\lambda +1) \dots (\lambda +n-1) $.
It appears not only in physics as special functions, but also in pure mathematics (e.g.~the number theory, the geometry including the conformal mapping). 
The Gauss hypergeometric function satisfies the hypergeometric differential equation 
\begin{align}
&  z(1-z) \frac{d^2y}{dz^2} + \left( \gamma - (\alpha + \beta +1)z \right) \frac{dy}{dz} -\alpha \beta  y=0,
\label{eq:GaussHGE}
\end{align}
which is a standard form of second order Fuchsian differential equation with three regular singularities $\{ 0,1,\infty \}$.

A $q$-difference analogue of the hypergeometric function was introduced as
\begin{align}
& _2 \phi _1 (a ,b ;c ;x) = \sum_{n=0}^{\infty} \frac{(a ;q)_n (b ;q)_n }{(q;q) _n (c ;q)_n } x^n, 
\label{eq:qhypser}
\end{align}
where $(\lambda ,q)_n $ is the $q$-Pochhammer symbol defined by $(\lambda ,q)_0 =1 $ and 
\begin{equation}
(\lambda ,q)_n= (1- \lambda )(1- \lambda q) (1- \lambda q^2) \dots (1- \lambda q^{n-1})  
\end{equation}
for the positive integer $n$.
Eq.~(\ref{eq:qhypser}) is called the basic hypergeometric function or the $q$-hypergeometric function.
It has been studied extensively from the 19th century, and a large number of formulas have been found.
The basic hypergeometric function $_2 \phi _1 (a ,b ;c ;x) $ satisfies the basic (or the $q$-difference) hypergeometric equation 
\begin{equation}
(x-q) f(x/q) - ((a+b)x -q-c)f(x)+ (abx-c)f(q x)=0. 
\label{eq:qhyp}
\end{equation}
Note that every coefficient of $ f(x/q)$, $f(x)$ and $f(q x) $ is linear in $x$.
It is known that the functions
\begin{align}
& x^{-\alpha } \! _2 \phi _1 (a ,aq/c ;aq/b ; cq/(abx)), \;  x^{-\beta } \! _2 \phi _1 (b ,bq/c ;bq/a ; cq/(abx))
\label{eq:qHGSinfty}
\end{align}
are solutions to Eq.~(\ref{eq:qhyp}), where $\alpha $ and $\beta $ satisfy $a=q^{\alpha }$ and $b=q^{\beta }$.
Hahn found other solutions in \cite{Hahn49} and the function  
\begin{align}
& \sum_{n=0}^{\infty} (abx/c ;q)_n \frac{(a ;q)_n (b ;q)_n }{(abq/c;q) _n } q^n (= \! _3 \phi _2 (a ,b,abx/c ; abq/c,0 ; q) ) 
\label{eq:qHGSPoch}
\end{align}
is a solution to Eq.~(\ref{eq:qhyp}).
By a suitable limit as $q\to 1$ of the $q$-difference hypergeometric equation, we may obtain the hypergeometric differential equation (\ref{eq:GaussHGE}).
On the $q$-difference equation as Eq.~(\ref{eq:qhyp}) (or Eq.~(\ref{eq:axgbxgcxg})), the points $x=0$ and $x=\infty $ are special because they are the accumulation points by the $q$-shift $x \to qx \to q^2x \to \dots $ (or $x \to q^{-1}x \to q^{-2}x \to \dots $).
On the other hand, in the framework of the differential equation on the Riemann sphere $\C \cup \{ \infty \} $ we can transform two points $x=0, \infty $ to arbitrary distinct two points by the linear fractional transformation.

In this paper, we propose two variants of the $q$-hypergeometric equation where the point $x=0$ is essentially non-singular, and show that they have several explicit solutions.
Before introducing the variants of the $q$-hypergeometric equation, we recall the $q$-Heun equation.
Hahn \cite{Hahn} introduced a $q$-difference analogue of Heun's differential equation of the form
\begin{equation}
a(x) g(x/q) + b(x) g(x) + c(x) g(qx) =0 ,
\label{eq:axgbxgcxg}
\end{equation}
where $a(x)$, $b(x)$, $c(x)$ are the polynomials such that $\deg_x a(x)= \deg_x c(x)=2 $, $a(0) \neq 0 \neq c(0) $ and $\deg _x b(x) \leq 2$, and it was rediscovered in \cite{TakR} as the equation $ A^{\langle 4 \rangle} g(x)= E g(x) $ where $A^{\langle 4 \rangle} $ is the fourth degeneration of Ruijsenaars-van Diejen operator of one variable and $E$ is a complex number.
The equation $ A^{\langle 4 \rangle} g(x)= E g(x) $ is written as
\begin{align}
& (x-q^{h_1 +1/2} t_1) (x- q^{h_2 +1/2} t_2) g(x/q)  \label{eq:RuijD5-0} \\
&  + q^{\alpha _1 +\alpha _2} (x - q^{l_1-1/2}t_1 ) (x - q^{l_2 -1/2} t_2) g(qx) \nonumber \\
&  -\{ (q^{\alpha _1} +q^{\alpha _2} ) x^2 + E x + q^{(h_1 +h_2 + l_1 + l_2 +\alpha _1 +\alpha _2 )/2 } ( q^{\beta /2}+ q^{-\beta/2}) t_1 t_2 \} g(x) =0, \nonumber
\end{align}
and it gives a standard parameterization of the $q$-Heun equation.
Note that the Ruijsenaars-van Diejen system (\cite{vD0,RuiN}) is a difference (or a relativistic) analogue of quantum Inozemtsev system, and the (non-degenerate) Ruijsenaars-van Diejen operator is written in terms of the elliptic theta function.

The secondly degenerated Ruijsenaars-van Diejen operator $ A^{\langle 2 \rangle} $ and the thirdly degenerated Ruijsenaars-van Diejen operator $ A^{\langle 3 \rangle} $ were also obtained in \cite{TakR}, and they were investigated from the viewpoint of the $q$-analogue of the apparent singularity in \cite{TakqH}.
For example, the equation $A^{\langle 3 \rangle} g(x) =Eg(x) $ is characterized as Eq.~(\ref{eq:axgbxgcxg}) such that $a(x)$, $b(x)$ and $c(x)$ are polynomials, $\deg_x a(x)= \deg_x c(x)=3 $, $\deg _x b(x) \leq 3$, $a(0)\neq 0 \neq c(0)$ and the singularity $x=\infty $ is essentially non-singular (see \cite{TakqH} for more details and also for the equation $A^{\langle 2 \rangle} g(x) =Eg(x) $).

In this paper, we derive a variant of the $q$-hypergeometric equation by imposing the condition that the origin $x=0$ is essentially non-singular on the $q$-Heun equation in Eq.~(\ref{eq:RuijD5-0}).
Then the condition is written as $\beta = \pm 1$ and $ E=-q^{(\alpha _1 + \alpha _2 +h_1+h_2+l_1+l_2)/2} \{ (q^{- h_2 }+q^{-l_2 })t_1 + (q^{- h_1 }+ q^{- l_1 }) t_2 \}$, which we explain in section \ref{sec:sing}.
Thus we introduce a variant of the $q$-hypergeometric equation by 
\begin{align}
& (x-q^{h_1 +1/2} t_1) (x - q^{h_2 +1/2} t_2) g(x/q) \label{eq:qhypervar1} \\
& + q^{\alpha _1 +\alpha _2} (x - q^{l_1-1/2}t_1 ) (x - q^{l_2 -1/2} t_2) g(q x) \nonumber  \\
&  -[ (q^{\alpha _1} +q^{\alpha _2} ) x^2 +E x + p ( q^{1/2}+ q^{-1/2}) t_1 t_2 ] g(x) =0, \nonumber \\
& p= q^{(h_1 +h_2 + l_1 + l_2 +\alpha _1 +\alpha _2 )/2 } , \quad E= -p \{ (q^{- h_2 }+q^{-l_2 })t_1 + (q^{- h_1 }+ q^{- l_1 }) t_2 \} \nonumber .
\end{align}
We call Eq.~(\ref{eq:qhypervar1}) the variant of the $q$-hypergeometric equation of degree two.
Then it is shown that there exist two dimensional convergent solutions to Eq.~(\ref{eq:qhypervar1}) whose basis are written as
\begin{equation*}
x^{\lambda } (1+ \tilde{c}_2 x^2+  \tilde{c}_3 x^3 + \dots ), \; x^{\lambda +1} (1 + \tilde{c}'_1 x+ \tilde{c}'_2 x^2 + \dots ),
\end{equation*}
where 
\begin{equation*}
\lambda = (h_1 +h_2 -l_1-l_2 -\alpha _1-\alpha _2+1 )/2 .
\end{equation*}
By the transformation $h(x)= x^{-\lambda } g(x)$, the $q$-difference equation for the function $h(x)$ has two-dimensional holomorphic solutions about $x=0$, and we may regard the point $x=0$ as non-singularity.

We investigate several explicit solutions to the variant of the $q$-hypergeometric equation of degree two which was given in Eq.~(\ref{eq:qhypervar1}).
To begin with, we look for the solution about $x=\infty $ written as
\begin{equation}
x^{-\rho ' } (1+ c_1 x^{-1} +c_2 x^{-2 }+ \dots +c_n x^{-n} + \dots ) ,
\label{eq:expandmu}
\end{equation}
which is analogous to Eq.~(\ref{eq:qHGSinfty}).
Then we can see directly that $\rho ' = \alpha _1$ or $\alpha _2$ by substituting Eq.~(\ref{eq:expandmu}) into Eq.~(\ref{eq:qhypervar1}).
Moreover we may calculate the coefficients $c_n$ and we find that the solution for the case $\rho ' = \alpha _1$ is written as
\begin{align*}
g_1 (x) & = x^{-\alpha _1 } \sum _{n=0}^{\infty} (q^{1/2} x^{-1})^n \frac{ ( q^{\lambda +\alpha _1 } ; q )_n }{(q^{\alpha _1  -\alpha _2 +1 } ; q )_n } \\
& \qquad \cdot \sum _{k=0}^n  \frac{(q^{ \lambda +\alpha _1  -h_2 +l_2 }; q)_k (q^{  \lambda +\alpha _1  -h_1 +l_1 } ;q)_{n-k }}{(q;q )_k (q;q)_{n-k}}  (q^{l_1} t_1)^k (q ^{l_2}t_2 ) ^{n-k} \nonumber \\
 & = x^{-\alpha _1 } \sum _{\ell =0}^{\infty} \sum _{k=0}^{\infty } \frac{( q^{\lambda + \alpha _1 } ; q )_{k+\ell } }{(q^{ \alpha _1  - \alpha _2 +1 } ; q )_{k+\ell } } \frac{(q^{ \lambda +\alpha _1  -h_2 +l_2 }; q)_k (q^{  \lambda +\alpha _1  -h_1 +l_1 } ;q)_{\ell }}{(q;q )_k (q;q)_{\ell }} \nonumber \\
& \qquad \cdot (q^{l_1 +1/2} t_1 x^{-1} )^k (q ^{l_2 +1/2 }t_2 x^{-1} ) ^{\ell } . \nonumber
\end{align*}
It is a specialization of the $q$-Appell function
\begin{equation*}
\Phi ^{(1)}(a;b,b';c;q;x_1,x_2)= \sum _{\ell = 0}^{\infty } \sum _{ k = 0}^{\infty } \frac{(a;q)_{\ell + k} (b;q)_{\ell } (b';q )_k }{(q;q)_{\ell } (q;q)_k (c;q)_{\ell + k}}x_1^{\ell} x_2^k .
\end{equation*}
We can show that the variant of the $q$-hypergeometric equation of degree two is obtained by considering a restriction of the $q$-difference equations which the $q$-Appell function satisfies.
For details, see section \ref{sec:sol} and the appendix.

We investigate other solutions to the variant of the $q$-hypergeometric equation of degree two.
To find the solution in a similar form to Eq.~(\ref{eq:qHGSPoch}), we set 
\begin{align*}
& g_2(x)= x^{\lambda } \sum _{n=0}^{\infty} c_n \Big( \frac{x}{q^{l_1-1/2} t_1} ;q \Big) _n  , \quad (c_0=1),
\end{align*}
and substitute it into Eq.~(\ref{eq:qhypervar1}).
Then we find that the coefficients are determined as
\begin{align*}
&  c_n = q^n \frac{(q^{\lambda +\alpha _1 };q )_n (q^{\lambda +\alpha _2 };q )_n }{(q^{h_1 - l_1 +1};q )_n (q^{h_2 - l_1 +1} t_2/t_1 ;q )_n (q;q)_n } .
\end{align*}
We show in Theorem \ref{thm:thm2} that the function 
\begin{align*}
& g_2 (x)= x^{\lambda } \: _3 \phi _2 \Big( \frac{x}{q^{l_1-1/2} t_1},q^{\lambda +\alpha _1 },q^{\lambda +\alpha _2 }; q^{h_1 - l_1 +1},q^{h_2 - l_1 +1} \frac{t_2}{t_1} ;q;q  \Big),
\end{align*}
is a solution to Eq.~(\ref{eq:qhypervar1}).
Note that it would be related to the big $q$-Jacobi polynomial (see \cite{KLS}).

We also investigate a solution to the variant of the $q$-hypergeometric equation of degree two in the form
\begin{align*}
& g_3 (x)= x^{-\alpha _1 } \sum _{n=0}^{\infty} c_n \Big( \frac{q^{h_1+1/2} t_1 }{x}  ;q \Big)_n , \quad (c_0=1).
\end{align*}
Then the coefficients are determined as
\begin{align*}
& c_n =q^n \frac{(q^{\lambda +\alpha _1 };q )_n }{(q^{h_1 - l_2 +1} t_1/t_2 ;q )_n  } \sum _{k=0}^n   \frac{(q^{\lambda -h_2+l_2+\alpha _1} ;q)_k q^{k(k+1)/2} (-q^{h_1-l_2}t_1/t_2)^k }{(q^{h_1-l_1+1};q)_{k} (q;q)_k (q;q)_{n-k} } .
\end{align*}
Thus we can show that the function 
\begin{align*}
& g_3(x)= x^{-\alpha _1 } \Big[ \sum _{n=0}^{\infty} \Big( \frac{q^{h_1+1/2} t_1 }{x}  ;q \Big)_n \frac{q^n (q^{\lambda +\alpha _1 };q )_n }{(q^{h_1 - l_2 +1} t_1/t_2 ;q )_n  } \\
& \qquad \qquad \qquad \qquad \cdot \sum _{k=0}^n   \frac{(q^{\lambda -h_2+l_2+\alpha _1} ;q)_k q^{k(k+1)/2} (-q^{h_1-l_2}t_1/t_2)^k }{(q^{h_1-l_1+1};q)_{k} (q;q)_k (q;q)_{n-k} }  \Big] 
\nonumber
\end{align*}
is a solution to Eq.~(\ref{eq:qhypervar1}), which we show in Theorem \ref{thm:thm3}.

We introduce another variant of the $q$-hypergeometric equation by
\begin{align}
& (x-q^{h_1 +1/2} t_1) (x- q^{h_2 +1/2} t_2) (x- q^{h_3 +1/2} t_3) g(x/q) \label{eq:var3qhyp0} \\
&  + q^{2\alpha  +1} (x - q^{l_1-1/2}t_1 ) (x - q^{l_2 -1/2} t_2) (x - q^{l_3 -1/2} t_3) g(qx) \nonumber \\
& + q^{\alpha } [ - (q + 1 ) x^3 + q^{1/2} \{ (q^{h_1} + q^{l_1 })t_1 + (q^{h_2} + q^{l_2 })t_2 + (q^{h_3} + q^{l_3 })t_3 \} x^2 \nonumber \\
& \quad - q^{(h_1+h_2+h_3+l_1+l_2+l_3 +1)/2} \{ (q^{- h_1 }+q^{-l_1 })t_2 t_3 + (q^{- h_2 }+ q^{- l_2 }) t_1 t_3\nonumber \\
& \quad  + (q^{- h_3 }+ q^{- l_3 }) t_1 t_2  \} x + q^{(h_1 +h_2 + h_3 + l_1 + l_2 + l_3 )/2 } ( q + 1 ) t_1 t_2t_3 ] g(x) =0, \nonumber
\end{align}
which we call the variant of the $q$-hypergeometric equation of degree three.
See section \ref{sec:sing} for its derivation.
The parameter $\alpha $ in Eq.~(\ref{eq:var3qhyp0}) can be eliminated by setting $g(x) = x ^{-\alpha } \tilde{g}(x) $.
Note that Eq.~(\ref{eq:var3qhyp0}) with the condition $\alpha =-1/2$ is realized by a specialization of a variant of the $q$-Heun equation $A^{\langle 3 \rangle} g(x) =Eg(x) $, where $ A^{\langle 3 \rangle} $ is the thirdly degenerated Ruijsenaars-van Diejen operator.

We give a conjecture of solutions to the variant of the $q$-hypergeometric equation of degree three, which suggest the hypergeometric structure of Eq.~(\ref{eq:var3qhyp0}). 
\begin{conj} \label{conj:var3}
Set $\nu = (h_1 +h_2 +h_3 -l_1-l_2-l_3 +1)/2$.
Let $(i,i',i'')$ be a permutation of $(1,2,3)$.\\
(i) The function 
\begin{align*}
& h(x)= x^{\nu -\alpha } \sum _{n=0}^{\infty} c_n \Big( \frac{x}{q^{l_i-1/2} t_i} ;q  \Big)_n  ,
\end{align*}
whose coefficients are determined by
\begin{align*}
& c_n =\frac{q^n  (q^{\nu };q )_n }{(q^{h_{i'} - l_i +1} t_{i'}/t_i ;q )_n (q^{h_{i''} - l_i +1} t_{i''}/t_i ;q )_n  } \\
& \qquad \cdot \sum _{k' =0}^n \sum _{k'' =0}^{n-k'} c_{n,k',k''} \Big( \frac{ - q^{h_{i'}-l_i}t_{i'}}{t_i} \Big)^{k'} \Big( \frac{- q^{h_{i''}-l_i}t_{i''}}{t_i} \Big)^{k''} , \nonumber \\
& c_{n,k',k''} =  q^{(k'+k'')(k' +k'' +1)/2} \frac{(q^{\nu -h_{i'}+l_{i'} } ;q)_{k'} (q^{\nu -h_{i''}+l_{i''} } ;q)_{k''} }{(q;q)_{n-k'-k''} (q;q)_{k'}(q;q)_{k''} (q^{h_i-l_i +1};q)_{k' + k''} } , \nonumber 
\end{align*}
is a solution to Eq.~(\ref{eq:var3qhyp0}).\\
(ii) The function 
\begin{align*}
& h(x)= x^{-\alpha } \sum _{n=0}^{\infty} c_n \Big( \frac{q^{h_i+1/2} t_i }{x}  ;q \Big)_n  , 
\end{align*}
whose coefficients are determined by
\begin{align*}
& c_n =\frac{q^n (q^{\nu };q )_n }{(q^{h_i - l_{i'} +1} t_i /t_{i'} ;q )_n (q^{h_i - l_{i''} +1} t_i/t_{i''} ;q )_n  } \\
& \qquad \cdot \sum _{k' =0}^n \sum _{k'' =0}^{n-k'} c_{n,k',k''} \Big( \frac{- q^{h_i-l_{i'}}t_i}{t_{i'}} \Big)^{k'} \Big( \frac{- q^{h_i-l_{i''}}t_i}{t_{i''}} \Big)^{k''} ,\\
& c_{n,k',k''} =  q^{(k'+k'')(k' +k'' +1)/2} \frac{(q^{\nu -h_{i'}+l_{i'} } ;q)_{k'} (q^{\nu -h_{i''}+l_{i''} } ;q)_{k''} }{(q;q)_{n-k'-k''} (q;q)_{k'}(q;q)_{k''} (q^{h_i-l_i +1};q)_{k' + k''} } , \nonumber 
\end{align*}
is a solution to Eq.~(\ref{eq:var3qhyp0}).
\end{conj}

This paper is organized as follows.
In section \ref{sec:sing}, we consider apparency of the singularity and two variants of the $q$-hypergeometric equation.
In section \ref{sec:sol}, we establish three solutions to the variant of the $q$-hypergeometric equation of degree two.
In section \ref{sec:limit}, we investigate the relationship between the $q$-hypergeometric equation and its variants, and consider the limit to obtain the differential equations as $q\to 1$.
In section \ref{sec:concl}, we give concluding remarks.
In the appendix, we investigate the $q$-Appell function.

Note that we treat the infinite sums formally and we do not discuss convergence.

\section{Apparency of the singularity and variants of the $q$-hypergeometric equation} \label{sec:sing}
 
We recall some facts on the regular singularity of the linear differential equation
\begin{equation}
y''  +p (z) y' +q (z) y=0 ,
\label{eq:y2}
\end{equation}
on the Riemann sphere ${\mathbb C} \cup \{ \infty \}$.
The point $z=0$ is called the regular singularity, if the functions $z p(z) $ and $z^2 q(z)$ are analytic about the point $z=0$.
Write $p_0= \lim _{z \to 0} z p(z)$ and $q_0 = \lim _{z \to 0} z^2 q(z)$.
Then the quadratic polynomial $\rho ^2 +(p_0 -1) \rho +q_0 $ and its roots are called the characteristic polynomial and the exponents about the regular singularity $z=0$.
They play important roles in the study of the local solutions about the regular singularity.
We assume that the origin $z=0$ is a regular singularity of the differential equation (\ref{eq:y2}) and write
\begin{equation*}
z p(z)= \sum _{m=0}^{\infty} p_m z^m, \; z^2 q (z)= \sum _{m=0}^{\infty} q_m z^m.
\end{equation*}
Let $\lambda _1$ and $\lambda _2$ be the exponents and write the characteristic polynomial as $F(\rho )$.
Then we have $F(\rho ) = \rho ^2 +(p_0 -1) \rho +q_0 = (\rho -\lambda _1 ) (\rho -\lambda _2 )$.
We now investigate the solution to the differential equation (\ref{eq:y2}) in the form
\begin{equation}
y= z^{\rho  }\sum _{m=0}^{\infty}  c_m z^m  , \quad c_0 =1 .
\label{eq:y2yseries}
\end{equation}
By substituting it into the differential equation, we find that $ F(\rho ) =0 $ and $F(\rho +m) c_m +  \sum _{k=1}^{m} \{ (\rho + m-k )  p_k +q_k \}  c_{m-k}=0 $ for $m=1,2,\dots $.
Therefore the value $\rho $ coincides with one of the exponents, and write $\rho  = \lambda _j$ for $j=1$ or $2$.
If $F(\lambda _j +m ) \neq 0$ for every positive integer $m$, then the coefficients $c_m $ $(m=1,2,\dots )$ are determined recursively and we obtain a convergent solution written as Eq.~(\ref{eq:y2yseries}).
If the difference of the exponents is not an integer (i.e. $\lambda _2- \lambda _1 \not \in {\mathbb Z}$), then we have two independent solutions.
On the other hand, if $\lambda _2- \lambda _1 \in {\mathbb Z}_{>0}$, then we have $F(\lambda _1 + N ) = 0$ for $N=\lambda _2 - \lambda _1 (>0)$ and the value $c_N$ is not determined by the recursive relation.
The values $c_m$ for $m\leq N-1$ are determined by the recursive relation and we have a necessary condition 
\begin{equation}
\sum _{k=1}^{N} \{ (\lambda _1 + N-k )  p_k +q_k \}  c _{N-k} = 0 
\label{eq:diffalappa}
\end{equation}
for the existence of the non-zero series solution corresponding to the exponent $\lambda _1$.
Otherwise, we need the logarithmic term for expressing the solution corresponding to the exponent $\lambda _1$.
The regular singularity $z=0$ is called apparent, if the condition in Eq.~(\ref{eq:diffalappa}) is satisfied.

We restrict to the case that $N=1$.
Then the condition for apparency of the singularity is written as $\lambda _1 p_1 +q_1 =0$.
It follows from $\lambda _2 - \lambda _1 =1$ that $p_0= 1-(\lambda _1 + \lambda _2  ) = -2\lambda _1$ and $q_0= \lambda _1 (\lambda _1 + 1)$.
Set $ u =z^{\lambda _1} y$.
Then the differential equation for $u$ is written as 
\begin{equation*}
\frac{d^2u}{dz^2} + \left\{ \frac{2\lambda _1}{z} +p(z)  \right\} \frac{du}{dz} +\left\{ \frac{\lambda _1 (\lambda _1 -1)}{z^2} + \frac{\lambda _1}{z} p(z) +q(z)  \right\} u=0. 
\end{equation*}
By recalling Eq.~(\ref{eq:diffalappa}), we see that the coefficients of $du/dz $ and $u$ are analytic about $z=0$.
Therefore we obtain the following proposition;
\begin{prop}
We assume that the exponents of the regular singularity $z=0$ on Eq.~(\ref{eq:y2}) are $\lambda _1$ and $\lambda _1 +1$, and the point $z=0$ is apparent singularity.
Set $u= z^{\lambda _1} y$.
Then the point $z=0$ is not the singularity on the differential equation for the function $u$.
\end{prop}
This motivates us to introduce the non-singularity for the origin $z=0$ of the difference equation.

We now investigate the linear difference equation
\begin{equation}
u(x) g(x/q) + v(x) g(x) + w(x) g(qx) =0 ,
\label{eq:axgbxgcxg1}
\end{equation}
where $u(x)$, $v(x)$, $w(x)$ are polynomials. 
Write 
\begin{align*}
& u(x) = \sum_{k=0}^{\infty } u_k x^k , \; v(x) = \sum_{k=0}^{\infty } v_k x^k , w(x) = \sum_{k=0}^{\infty } w_k x^k ,
\end{align*}
$n_k=v_k=w_k=0 $ for sufficiently large $k$ and assume that $u_0 \neq 0 \neq w_0$.
Then the origin $x=0$ may be regarded as the regular singularity (see \cite{Adams28,TakqH}).
We investigate local solutions about $x=0$ of the form
\begin{align*}
& g(x)= x^{\rho } \sum _{k =0}^{\infty} c_{k } x^{k }, \; c_0 =1.
\end{align*}
By substituting this into Eq.~(\ref{eq:axgbxgcxg1}) and looking at the coefficient of $x^{\rho  }$, we have
\begin{align*}
& q^{-\rho } u_0 + v_0 + q^{\rho } w_0 =0,
\end{align*}
and we call it the characteristic equation about $x=0$.
The exponents are the values $\rho $ which satisfy the characteristic equation.
Let $\lambda '$ be an exponent.
If $\lambda ' +n$ ($n =1,2,\dots  $) is not an exponent, then the coefficient $c_n$ ($n =1,2,\dots  $) is determined by 
\begin{align*}
&  ( q^{-\lambda ' -n} u_0 + v_0 + q^{\lambda '+n} w_0  ) c_{n} = - \sum_{k =0}^{n-1} ( q^{-\lambda '-k } u_{n-k } + v_{n-k } + q^{\lambda '+k } w_{n-k } ) c_{k } 
\end{align*}
(see \cite{TakqH}).
If $\lambda ' +N$ is an exponent for some positive integer $N$, then we need the equation 
\begin{align}
& \sum_{k =0}^{N-1} ( q^{-\lambda ' -k } u_{N-k } + v_{N-k } + q^{\lambda ' +k } w_{N-k } ) c_{k } =0 .
\label{eq:appsing}
\end{align}
to have series solutions, otherwise we need logarithmic terms for the solution.
If Eq.~(\ref{eq:appsing}) is satisfied, then the singularity $x=0$ is called apparent (or non-logarithmic).

We restrict to the case $N=1$.
Then the value $\lambda '+1 $ is also the exponent, i.e. we have
\begin{align*}
& q^{-\lambda '} u_0 + v_0 + q^{\lambda '} w_0 =0, \;  q^{-\lambda ' -1} u_0 + v_0 + q^{\lambda ' +1} w_0 =0.
\end{align*}
In this case, Eq.~(\ref{eq:appsing}) is written as
\begin{align}
& q^{-\lambda '  } u_{1 } + v_{1 } + q^{\lambda ' } w_{1 }  =0 .
\label{eq:appsing1}
\end{align}
We use this condition for criterion of elimination of the singularity about $x=0$.

We examine it for the $q$-Heun equation.
Recall that the $q$-Heun equation is written as
\begin{align*}
& (x-q^{h_1 +1/2} t_1) (x- q^{h_2 +1/2} t_2) g(x/q)  \\
& + q^{\alpha _1 +\alpha _2} (x - q^{l_1-1/2}t_1 ) (x - q^{l_2 -1/2} t_2) g(qx) \nonumber \\
& -\{ (q^{\alpha _1} +q^{\alpha _2} ) x^2 + E x + q^{(h_1 +h_2 + l_1 + l_2 +\alpha _1 +\alpha _2 )/2 } ( q^{\beta /2}+ q^{-\beta/2}) t_1 t_2 \} g(x) =0. \nonumber
\end{align*}
The characteristic equation about the origin $x=0$ is written as
\begin{align*}
q^{h_1 + h_2 + 1} q^{-\lambda } -q^{(h_1 +h_2 + l_1 + l_2 +\alpha _1 +\alpha _2 )/2 } ( q^{\beta /2}+ q^{-\beta/2})+ q^{l_1 +l_2 +\alpha _1 +\alpha _2 -1 } q^{\lambda } =0. 
\end{align*}
Hence the values 
\begin{equation*}
(h_1 +h_2 -l_1-l_2 -\alpha _1-\alpha _2 -\beta +2)/2, \; (h_1 +h_2 -l_1-l_2 -\alpha _1-\alpha _2 +\beta +2)/2 
\end{equation*}
are exponents about $x=0$.
We assume that the values $\lambda $ and $\lambda +1$ are exponents about $x=0$ and the singularity $x=0$ is apparent.
Then we have $\beta =\pm 1$, $\lambda = (h_1 +h_2 -l_1-l_2 -\alpha _1-\alpha _2+1 )/2 $ and the condition in Eq.~(\ref{eq:appsing1}), which is equivalent to 
\begin{equation*}
 E=-q^{(\alpha _1 + \alpha _2 +h_1+h_2+l_1+l_2)/2} \{ (q^{- h_2 }+q^{-l_2 })t_1 + (q^{- h_1 }+ q^{- l_1 }) t_2 \} ,
\end{equation*}
is satisfied.
Thus we obtain the variant of the $q$-hypergeometric equation of degree two, which was introduced in Eq.~(\ref{eq:qhypervar1}).

Next we examine the $q$-difference equation of the form
\begin{align*}
& (x-q^{h_1 +1/2} t_1) (x- q^{h_2 +1/2} t_2) (x- q^{h_3 +1/2} t_3) g(x/q) \\
&  + \tilde{b} (x - q^{l_1-1/2}t_1 ) (x - q^{l_2 -1/2} t_2) (x - q^{l_3 -1/2} t_3) g(qx) \nonumber \\
& - [ b_3 x^3 + b_2 x^2 +b_1 x +b_0 ] g(x)=0 .\nonumber 
\end{align*}
We investigate the condition that the difference of the exponents at $x=0$ is one and the singularity $x=0$ is apparent and also a similar condition for $x=\infty $.

We investigate local solutions about $x= \infty$ in the form
\begin{align*}
& h(x)= (1/x)^{\tilde{\rho }} \sum _{k =0}^{\infty} \tilde{c}_{k } (1/x)^{k }, \; \tilde{c}_0 \neq 0.
\end{align*}
Then a necessary condition for existence of non-zero solution is given by
\begin{align*}
& q^{\tilde{\rho }} -b_3 + \tilde{b} q^{-\tilde{\rho }}  =0,
\end{align*}
which is the characteristic equation about $x=\infty $ and we call the roots the exponents (see \cite{TakqH}).
We assume that $\alpha $ and $\alpha  +1 $ are exponents and the singularity $x =\infty $ is apparent.
Then the condition is written as
\begin{align*}
& q^{\alpha } -b_3 + \tilde{b} q^{-\alpha  }  =0, \; q^{\alpha  +1} -b_3 + \tilde{b} q^{-\alpha  -1}  =0, \\
& ( q^{h_1 +1/2} t_1 + q^{h_2 +1/2} t_2 + q^{h_3 +1/2} t_3) q^{\alpha  } + b _2 \nonumber \\
& \qquad + \tilde{b} ( q^{l_1-1/2}t_1 + q^{l_2 -1/2} t_2 + q^{l_3 -1/2} t_3)q^{-\alpha  } =0 . \nonumber 
\end{align*}
Hence we have $\tilde{b} = q^{2\alpha  +1} $, $b_3 =  q^{\alpha  } + q^{\alpha  +1}  $ and $b_2 = -( q^{h_1 } t_1 + q^{h_2 } t_2 + q^{h_3 } t_3 + q^{l_1 }t_1 + q^{l_2 } t_2 + q^{l_3 } t_3) q^{\alpha +1/2}$.
Next we consider the origin $x=0$.
The characteristic equation about $x= 0$ is written as 
\begin{align*}
& - q^{-\rho  } q^{h_1 +h_2 +h_3 + 3/2} t_1 t_2 t_3 - b_0 - q^{2\alpha  +1} q^{\rho } q^{l_1 +l_2 +l_3 - 3/2} t_1 t_2 t_3 =0.  
\end{align*}
The condition that $\lambda $ and $\lambda  +1 $ are exponents and the singularity $x =0 $ is apparent is written as
\begin{align*}
& q^{-\lambda } q^{h_1 +h_2 +h_3 + 3/2} t_1 t_2 t_3 +b_0 + q^{2\alpha  +1} q^{\lambda } q^{l_1 +l_2 +l_3 - 3/2} t_1 t_2 t_3 =0, \\
& q^{-\lambda -1} q^{h_1 +h_2 +h_3 + 3/2} t_1 t_2 t_3 +b_0 + q^{2\alpha  +1} q^{\lambda +1 } q^{l_1 +l_2 +l_3 - 3/2} t_1 t_2 t_3 =0, \nonumber \\
& q^{-\lambda } q^{h_1 +h_2 +h_3 + 3/2}t_1 t_2 t_3 (1/(q^{h_1 +1/2} t_1) + 1/(q^{h_2 +1/2} t_2) + 1/(q^{h_3 +1/2} t_3) ) - b_1 \nonumber \\
& + q^{2\alpha  +1} q^{\lambda } q^{l_1 +l_2 +l_3 - 3/2} t_1 t_2 t_3  (1/(q^{l_1 -1/2} t_1) + 1/(q^{l_2 -1/2} t_2) + 1/(q^{l_3 -1/2} t_3) ) =0 . \nonumber
\end{align*}
Thus we have $\lambda = (h_1 +h_2 +h_3 -l_1 -l_2 -l_3 - 2\alpha + 1)/2 $, $b_0 = -  q^{(h_1 +h_2 +h_3 +l_1 +l_2 +l_3 + 2\alpha )/2 } (1 +q) t_1 t_2 t_3 $, $b_1 = q^{(h_1 +h_2 +h_3 +l_1 +l_2 +l_3 + 2\alpha +1 )/2 } t_1 t_2 t_3 (1/(q^{h_1 } t_1) + 1/(q^{h_2 } t_2) + 1/(q^{h_3 } t_3) + 1/(q^{l_1 } t_1) + 1/(q^{l_2 } t_2) + 1/(q^{l_3 } t_3) )$.
Therefore we obtain the variant of the $q$-hypergeometric equation of degree three, which was introduced in Eq.~(\ref{eq:var3qhyp0}).
Set $\nu = (h_1 +h_2 +h_3 -l_1-l_2-l_3 +1)/2$.
The exponents about $x=0$ are $\nu -\alpha $ and $ \nu -\alpha +1$, and the singularity $x=0$ is apparent.
The exponents about $x=\infty $ are $\alpha $ and $\alpha +1$, and the singularity $x=\infty $ is apparent.

\section{Three solutions to the variant of the $q$-hypergeometric equation of degree two} \label{sec:sol}

In this section, we obtain three solutions to the variant of the $q$-hypergeometric equation of degree two.
One of them is related with the $q$-Appell function defined in the appendix.
In Proposition \ref{prop:qApp2ndorder}, it is shown that the specialized $q$-Appell function
\begin{equation}
f(x_1,x_2 )=\sum_{m=0}^\infty \sum_{n=0}^\infty \frac{(a;q)_{m+n}(b;q)_m (b';q)_n}{(b b' ;q)_{m+n}(q;q)_m(q;q)_n}(x_1)^m (x_2)^n
\label{eq:q-Appbb'c}
\end{equation}
satisfies the $q$-difference equation
\begin{align}
& (aq x_1-b')(aq x_2-b)f(q^2 x_1,q^2 x_2 ) +q(x_1-1)(x_2 -1)f(x_1,x_2)  \label{eq:2q-Appell} \\
& \qquad -(aq(q+1)x_1 x_2 -q(a+b)x_1 -q(a+b')x_2 +b b' +q)f(qx_1,qx_2 ) = 0 . \nonumber 
\end{align}
We can obtain the variant of the $q$-hypergeometric equation of degree two (see Eq.~(\ref{eq:qhypervar1})) from Eq.~(\ref{eq:2q-Appell}) as follows;
\begin{prop} \label{prop:qAppeqn}
Set
\begin{align} 
& a=q^{\lambda + \alpha_1 }, \; b= q^{\lambda + \alpha_1 +l_2 -h_2 }, \; b'= q^{\lambda + \alpha_1 +l_1 -h_1 }, \label{eq:abb'la} \\
& \lambda = (h_1 +h_2 -l_1-l_2 -\alpha _1-\alpha _2+1 )/2 . \nonumber
\end{align}
If the function $f(x_1,x_2) $ satisfies Eq.~(\ref{eq:2q-Appell}), then the function 
\begin{equation*}
g(x) = x ^{-\alpha_1 } f( q^{l_1+ 1/2} t_1 /x , q^{l_2+ 1/2} t_2 /x ) 
\end{equation*}
satisfies Eq.~(\ref{eq:qhypervar1}).
\end{prop}
\begin{prf}
We substitute $x_1 =d_1 /x$, $x_2=d_2 /x$, $f(qx_1,qx_2)=x^{d_3} g(x)$ into Eq.~(\ref{eq:2q-Appell}). 
Then we have
\begin{align*}
&(x-ab'^{-1}qd_1)(x-ab^{-1}qd_2)g(x/q) \\
&-q^{d_3} (b b' ) ^{-1}\{(b b' +q)x^2-q(a+b)d_1 x -q(a+b')d_2 x +a (q^2+q)d_1d_2 \}g(x) \nonumber \\
&+q^{2d_3+1}(b b' )^{-1}(x -d_1)(x -d_2)g(q x) =0 . \nonumber 
\end{align*}
By setting $a=q^{\lambda + \alpha_1 }$, $b= q^{\lambda + \alpha_1 +l_2 -h_2 }$, $b'= q^{\lambda + \alpha_1 +l_1 -h_1 }$, $d_1= q^{l_1-1/2}t_1 $, $d_2=q^{l_2-1/2}t_2$ and $d_3=\alpha_1$, we obtain Eq.~(\ref{eq:qhypervar1}).
\end{prf}

A solution to the variant of the $q$-hypergeometric equation of degree two, which is related with the $q$-Appell function, is written as follows;
\begin{thm} \label{thm:thm1}
Let $\lambda = (h_1 +h_2 -l_1-l_2 -\alpha _1-\alpha _2+1 )/2$.
The function 
\begin{align}
g (x) = & x^{-\alpha _1 } \sum _{n=0}^{\infty} (q^{1/2} x^{-1})^n \frac{ ( q^{\lambda +\alpha _1 } ; q )_n }{(q^{\alpha _1  -\alpha _2 +1 } ; q )_n } \label{eq:g1x} \\
& \cdot \sum _{k=0}^n  \frac{(q^{ \lambda +\alpha _1  -h_2 +l_2 }; q)_k (q^{  \lambda +\alpha _1  -h_1 +l_1 } ;q)_{n-k }}{(q;q )_k (q;q)_{n-k}}  (q^{l_1} t_1)^k (q ^{l_2}t_2 ) ^{n-k} \nonumber 
\end{align}
is a solution of the variant of the $q$-hypergeometric equation of degree two (i.e.~Eq.~(\ref{eq:qhypervar1})).
\end{thm}
\begin{prf}
The function $f(x_1,x_2 )$ in Eq.~(\ref{eq:q-Appbb'c}) satisfies Eq.~(\ref{eq:2q-Appell}). 
It follows from Proposition \ref{prop:qAppeqn} that the function $g(x) = x ^{-\alpha_1 } f( q^{l_1+ 1/2} t_1 /x , q^{l_2+ 1/2} t_2 /x ) $ satisfies Eq.~(\ref{eq:qhypervar1}) by specializing the parameters as Eq.~(\ref{eq:abb'la}).
Then the function $g(x)$ is expressed as Eq.~(\ref{eq:g1x}) in this specialization, and we obtain the theorem.
\end{prf}
\\
Note that the function which is obtained by replacing $\alpha _1 $ with $ \alpha _2$ is also a solution of the variant of the $q$-hypergeometric equation of degree two.
\begin{thm} \label{thm:thm2}
Set $(i,i')=(1,2) $ or $(2,1)$.
Then the function
\begin{align*}
& g (x)= x^{\lambda } \sum _{n=0}^{\infty} ( x/(q^{l_i-1/2} t_i) ;q )_n \frac{(q^{\lambda +\alpha _1 };q )_n (q^{\lambda +\alpha _2 };q )_n }{(q^{h_i - l_i +1};q )_n (q^{h_{i'} - l_i +1} t_{i'} /t_i ;q )_n (q;q)_n } q^n  
\end{align*}
is a solution of the variant of the $q$-hypergeometric equation of degree two (i.e.~Eq.~(\ref{eq:qhypervar1})).
\end{thm}
\begin{prf}
We show the theorem in the case $(i,i')=(1,2)$.
Note that where the case  $(i,i')=(2,1)$ is similar.

We determine the solution of the variant of the $q$-hypergeometric equation of degree two in the form
\begin{equation}
g(x)=x^{\lambda}\sum^{\infty}_{n=0}a_{n}\left(\frac{x}{q^{l_1-1/2} t_1 };q\right)_{n}, \quad a_0=1.
\label{eq:gxg2}
\end{equation}
Set $a_{-1}= a_{-2} =a_{-3} =0$.
We substitute it into Eq.~(\ref{eq:qhypervar1}).
On the term containing $g(qx)$, we have
\begin{align}
&q^{\alpha _1 +\alpha _2}(x-q^{l_1-\frac{1}{2}}t_1)(x-q^{l_2-\frac{1}{2}}t_2)g(qx) \label{eq:gqx} \\
&=q^{\alpha_1+\alpha_2+l_1+l_2-1}t_1t_2\biggl(1-\frac{x}{q^{l_1-1/2}t_1}\biggr)\biggl(1-\frac{x}{q^{l_2-1/2}t_2}\biggr)(qx)^{\lambda }\sum^{\infty}_{n=0}a_{n}\left(\frac{qx}{q^{l_1-1/2} t_1 };q\right)_{n} \nonumber \\
&=q^{\alpha_1+\alpha_2+l_1+l_2-1}t_1t_2\biggl(1-\frac{x}{q^{l_2-1/2}t_2}\biggr)(qx)^{\lambda }\sum^{\infty}_{n=0}a_{n}\left(\frac{x}{q^{l_1-1/2} t_1 };q\right)_{n+1} . \nonumber 
\end{align}
It follows from the identity 
\begin{equation*}
\biggl(1-\frac{x}{q^{l_2-1/2}t_2}\biggr) =(1- q^{l_1-l_2-n-1} t_1/t_2 ) + q^{l_1-l_2-n-1} t_1/t_2 \biggl(1-\frac{q^{n+1}x}{q^{l_1-1/2}t_1}\biggr) 
\end{equation*}
that Eq.~(\ref{eq:gqx}) is equal to 
\begin{align*}
& x^{\lambda }q^{\alpha_1+\alpha_2+\lambda +l_1+l_2-1}t_1t_2\sum^{\infty}_{n=0}\biggl\{a_{n} (1- q^{l_1-l_2-n-1} t_1/t_2 )\left(\frac{x}{q^{l_1-1/2} t_1 };q\right)_{n+1} \\
& \qquad \qquad \qquad \qquad \qquad \qquad + a_{n} q^{l_1-l_2-n-1} t_1/t_2 \left(\frac{x}{q^{l_1-1/2} t_1 };q\right)_{n+2}\biggr\}  \nonumber \\
&=x^{\lambda }q^{\alpha_1+\alpha_2+\lambda +l_1+l_2-1}t_1t_2\sum^{\infty}_{n=0}\biggl\{ a_{n-1}(1- q^{l_1-l_2-n} t_1/t_2 )\left(\frac{x}{q^{l_1-1/2} t_1 };q\right)_{n}  \nonumber \\
& \qquad \qquad \qquad \qquad \qquad \qquad + a_{n-2} q^{l_1-l_2-n+1} t_1/t_2 \left(\frac{x}{q^{l_1-1/2} t_1 };q\right)_{n}\biggr\} .  \nonumber 
\end{align*}
On the term containing $g(x/q)$, we have
\begin{align}
&(x-q^{h_1 +\frac{1}{2}} t_1)(x-q^{h_2 +\frac{1}{2}}t_2)g(x/q) \label{eq:gx/q}   \\
&=q^{h_1+h_2+1 }t_1t_2\biggl(1-\frac{x}{q^{h_1+1/2}t_1}\biggr)\biggl(1-\frac{x}{q^{h_2+1/2}t_2}\biggr)\biggl(\frac{x}{q}\biggr)^{\lambda } \sum^{\infty}_{n=0}a_{n}\left(\frac{xq^{-1}}{q^{l_1-1/2}};q\right)_{n}  \nonumber \\
&=q^{h_1+h_2+1}t_1t_2 \biggl(\frac{x}{q}\biggr)^{\lambda } \biggl(1-\frac{x}{q^{h_1+1/2}t_1}\biggr)\biggl(1-\frac{x}{q^{h_2+1/2}t_2}\biggr) \nonumber \\
& \qquad \qquad \qquad \qquad \qquad \qquad \qquad  \cdot \biggl(1-\frac{xq^{-1}}{q^{l_1-1/2}t_1}\biggr) \sum^{\infty}_{n=0}a_{n}\left(\frac{x}{q^{l_1-1/2}};q\right)_{n-1} .  \nonumber 
\end{align}
It follows from the identity 
\begin{align*}
&\biggl(1-\frac{x}{q^{h_1+1/2}t_1}\biggr)\biggl(1-\frac{x}{q^{h_1+1/2}t_2}\biggr) \biggl(1-\frac{xq^{-1}}{q^{l_1-1/2}t_1}\biggr) \\
&=A^{[-]}_{n}+B^{[-]}_{n}\biggl(1-\frac{q^{n-1}x}{q^{l_1-1/2}t_1}\biggr)+C^{[-]}_{n}\biggl(1-\frac{q^{n-1}x}{q^{l_1-1/2}t_1}\biggr)\biggl(1-\frac{q^{n}x}{q^{l_1-1/2}t_1}\biggr)  \nonumber \\
&\qquad \qquad \qquad \qquad +D^{[-]}_{n}\biggl(1-\frac{q^{n-1}x}{q^{l_1-1/2}t_1}\biggr)\biggl(1-\frac{q^{n}x}{q^{l_1-1/2}t_1}\biggr)\biggl(1-\frac{q^{n+1}x}{q^{l_1-1/2}t_1}\biggr) ,  \nonumber 
\end{align*}
\begin{align*}
& A^{[-]}_{n}=(1-q^{l_1-h_1-n})(1-q^{l_1-h_2-n} t_1/t_2)(1-q^{-n}) , \\
& B^{[-]}_{n}=q^{l_1-h_1-n}+q^{l_1-h_2-n} t_1/t_2 +q^{-n} -\{ q^{l_1-h_2-2n} (1+q^{l_1-h_1}) t_1/t_2  \nonumber \\
& \qquad + q^{l_1-h_1-2n} \}(1+q^{-1}) + q^{2l_1-h_1-h_2+3n}(1+q^{-1}+q^{-2}) t_1/t_2 , \nonumber \\
& C^{[-]}_{n}=\{q^{l_1-h_1-2n}+q^{l_1-h_2-2n}(1+q^{l_1-h_1}) t_1/t_2 \}q^{-1} \nonumber \\
& \qquad - q^{2l_1-h_1-h_2 - 3n-1}(1+q^{-1}+q^{-2}) t_1/t_2, \nonumber \\
& D^{[-]}_{n}= q^{2l_1-h_1-h_2-3n-3} t_1/t_2 ,  \nonumber 
\end{align*}
that Eq.~(\ref{eq:gx/q}) is equal to 
\begin{align*}
& x^{\lambda }q^{h_1+h_2-\lambda +1}t_1t_2\biggl[\sum^{\infty}_{n=0}a_{n}\biggl\{A^{[-]}_{n}\left(\frac{x}{q^{l_1-1/2} t_1 };q\right)_{n-1}+B^{[-]}_{n} \left(\frac{x}{q^{l_1-1/2} t_1 };q\right)_{n} \\
&\qquad \qquad \qquad \qquad +C^{[-]}_{n}\left(\frac{x}{q^{l_1-1/2} t_1 };q\right)_{n+1}+D^{[-]}_{n}\left(\frac{x}{q^{l_1-1/2} t_1 };q\right)_{n+2}\biggr\}\biggr] \nonumber  \\
&=x^{\lambda }q^{h_1+h_2-\lambda +1}t_1t_2\bigg\{\sum^{\infty}_{n=0}a_{n+1}A^{[-]}_{n+1}\left(\frac{x}{q^{l_1-1/2} t_1 };q\right)_{n}+a_{n}B^{[-]}_{n} \left(\frac{x}{q^{l_1-1/2} t_1 };q\right)_{n}  \nonumber \\
&\qquad \qquad \qquad \qquad +a_{n-1}C^{[-]}_{n-1}\left(\frac{x}{q^{l_1-1/2}};q\right)_{n}+a_{n-2}D^{[-]}_{n-2}\left(\frac{x}{q^{l_1-1/2} t_1 };q\right)_{n}\biggr\} . \nonumber  
\end{align*}
To calculate the term containing $g(x)$, we set
\begin{align*}
&(q^{\alpha _1} +q^{\alpha_2})x^2 + Ex +q^{(h_1 +h_2 +l_1+l_2+\alpha _1 +\alpha _2)/ 2}(q^{1/2}+q^{-1/2})t_1 t_2 \\
&=A^{[0]}_{n}+B^{[0]}_{n}\biggl(1-\frac{q^{n}x}{q^{l_1-1/2}t_1}\biggr)
+C^{[0]}_{n}\biggl(1-\frac{q^{n+1}x}{q^{l_1-1/2}t_1}\biggr)\biggl(1-\frac{q^{n}x}{q^{l_1-1/2}t_1}\biggr) . \nonumber 
\end{align*}
Then we have
\begin{align*}
& A^{[0]}_{n}=t_1t_2q^{h_1+h_2-\lambda_1+1/2}(q^{1/2}+q^{-1/2})+t_1q^{l_1-n-1/2}E+t_{1}^{2}q^{2l_1-2n-1}(q^{\alpha_1}+q^{\alpha_2}) , \\
& B^{[0]}_{n}=-t_1q^{l_1-n-1/2}E-t_{1}^{2}q^{2l_1-2n-2}(q^{\alpha_1}+q^{\alpha_2})(1+q) , \nonumber  \\
& C^{[0]}_{n}=t_{1}^{2}q^{2l_1-2n-2}(q^{\alpha_1}+q^{\alpha_2}) . \nonumber 
\end{align*}
Note that $ E= -q^{(h_1 +h_2 + l_1 + l_2 +\alpha _1 +\alpha _2 )/2 } \{ (q^{- h_2 }+q^{-l_2 })t_1 + (q^{- h_1 }+ q^{- l_1 }) t_2 \}  $.
Hence
\begin{align*}
&\{(q^{\alpha _1} +q^{\alpha_2})x^2 + Ex +q^{(h_1 +h_2 +l_1+l_2+\alpha _1 +\alpha _2)/2}(q^{1/2}+q^{- 1/2})t_1 t_2\} \\
&\qquad \cdot x^{\lambda } \sum^{\infty}_{n=0}a_{n}\left(\frac{x}{q^{l_1-1/2} t_1 };q\right)_{n} = x^{\lambda }\biggl\{ \sum^{\infty}_{n=0}a_{n}A^{[0]}_{n}\left(\frac{x}{q^{l_1-1/2} t_1 };q\right)_{n}
 \nonumber \\
& \qquad +\sum^{\infty}_{n=0}a_{n-1}B^{[0]}_{n-1}\left(\frac{x}{q^{l_1-1/2} t_1 };q\right)_{n} +\sum^{\infty}_{n=0}a_{n-2}C^{[0]}_{n-2}\left(\frac{x}{q^{l_1-1/2} t_1 };q\right)_{n} \biggr\} . \nonumber 
\end{align*}
Therefore, if $g(x)$ is expressed as Eq.~(\ref{eq:gxg2}), then we have
\begin{align*}
& (x-q^{h_1 +1/2} t_1) (x - q^{h_2 +1/2} t_2) g(x/q)  \\
&  \quad + q^{\alpha _1 +\alpha _2} (x - q^{l_1-1/2}t_1 ) (x - q^{l_2 -1/2} t_2) g(q x) \nonumber \\
&  \quad -[ (q^{\alpha _1} +q^{\alpha _2} ) x^2 +E x + q^{(h_1 +h_2 + l_1 + l_2 +\alpha _1 +\alpha _2 )/2 } ( q^{1/2}+ q^{-1/2}) t_1 t_2 ] g(x) \nonumber \\
& = \sum _{n=0}^{\infty} -t_1^{2}q^{2l_1-\lambda -3n-2} Q(a_{n+1},a_n,a_{n-1},a_{n-2} ) \left(\frac{x}{q^{l_1-1/2} t_1 };q\right)_{n} , \nonumber 
\end{align*}
where
\begin{align*}
& Q(a_{n+1},a_n,a_{n-1},a_{n-2} ) =   \\
& \qquad \qquad \qquad (1-q^{h_1-l_1+n+1})(1-q^{h_2-l_1+n+1} t_2/t_1)(1-q^{n+1})a_{n+1} \nonumber \\
& \qquad \qquad \qquad -q(1-q^{\lambda +\alpha_1+n})(1-q^{\lambda +\alpha_2 + n})a_{n}  \nonumber \\
& \qquad \qquad \qquad -q^{3}(1+q^{-1})(1-q^{h_1-l_1+n})(1-q^{h_2-l_1+n}t_2/t_1)(1-q^{n})a_{n} \nonumber \\
& \qquad \qquad \qquad +q^{4}(1+q^{-1})(1-q^{\lambda +\alpha_1+n-1})(1-q^{\lambda +\alpha_2+n-1})a_{n-1} \nonumber \\
& \qquad \qquad \qquad +q^{5}(1-q^{h_1-l_1+n-1})(1- q^{h_2-l_1+n-1} t_2/t_1 )(1-q^{n-1})a_{n-1} \nonumber \\
& \qquad \qquad \qquad -q^{6}(1-q^{\lambda +\alpha_1+n-2})(1-q^{\lambda +\alpha_2+n-2})a_{n-2} . \nonumber 
\end{align*}
Then $Q(a_{0},0 ,0 ,0 )$ is identically zero, and the coefficients $a_{n+1}$ $(n=0,1,2,\dots )$ are determined recursively by $Q(a_{n+1},a_n,a_{n-1},a_{n-2} ) =0$.
We set
\begin{align*}
& \tilde{a} _{n}=(1-q^{h_1-l_1+n+1})(1- q^{h_2-l_1+n+1} t_2/t_1 )(1-q^{n+1})a_{n+1}\\
& \quad -q(1-q^{\lambda +\alpha_1+n}) (1-q^{\lambda +\alpha_2+n})a_{n} . \nonumber
\end{align*}
Then the equation $Q(a_{n+1},a_n,a_{n-1},a_{n-2} ) =0 $ is equivalent to 
\begin{equation*}
\tilde{a} _{n}-q^{2}(1+q)\tilde{a} _{n-1}+q^{5}\tilde{a} _{n-2}=0, \quad n=0,1,2,\dots . 
\end{equation*}
It follows from $a_{-1}=a_{-2} =a_{-3} =0$ that $\tilde{a} _{-1} = \tilde{a} _{-2}=0$.
Hence we have $\tilde{a} _n =0$ for $n=0,1,2,\dots $. Namely
\begin{equation*}
a_{n+1}=\frac{q(1-q^{\lambda +\alpha_1+n})(1-q^{\lambda +\alpha_2+n})}
{(1-q^{h_1-l_1+n+1})(1- q^{h_2-l_1+n+1} t_2/t_1 )(1-q^{n+1})}a_{n}.
\end{equation*}
Therefore we obtain the solution
\begin{equation*}
g(x)=x^{\lambda } \sum^{\infty}_{n=0} 
q^{n}\frac{(q^{\lambda +\alpha_1};q)_n(q^{\lambda +\alpha_2};q)_n}{(q^{h_1-l_1+1};q)_n(q^{h_2-l_1+1}t_2/t_1;q)_n(q;q)_n}
\left(\frac{x}{q^{l_1-1/2}t_1 };q \right)_n .
\end{equation*}
\end{prf}
\begin{thm} \label{thm:thm3}
Set $(i,i')=(1,2) $ or $(2,1)$.
Then the function
\begin{align*}
& g (x)= x^{-\alpha _1 } \Big[ \sum _{n=0}^{\infty} ( q^{h_i+1/2} t_i /x  ;q )_n \frac{(q^{\lambda +\alpha _1 };q )_n }{(q^{h_i - l_{i'} +1} t_i/t_{i'} ;q )_n  } q^n \\
& \qquad \qquad \qquad \qquad \cdot \sum _{k=0}^n   \frac{(q^{\lambda -h_{i'} +l_{i'} +\alpha _1} ;q)_k }{(q^{h_i-l_i+1};q)_{k} (q;q)_k (q;q)_{n-k} } q^{k(k+1)/2} (-q^{h_i-l_{i'} }t_i/t_ {i'})^k  \Big] 
\nonumber
\end{align*}
is a solution of the variant of the $q$-hypergeometric equation of degree two (i.e.~Eq.~(\ref{eq:qhypervar1})).
\end{thm}
Note that the function which is obtained by replacing $\alpha _1 $ with $ \alpha _2$ is also a solution of the variant of the $q$-hypergeometric equation of degree two.
\\
\begin{prf}
We show the theorem in the case $(i,i')=(1,2)$.
We determine the solution of the variant of the $q$-hypergeometric equation of degree two in the form
\begin{equation}
g(x)=x^{-\alpha_1}\sum^{\infty}_{n=0}c_{n}\left(\frac{q^{h_1+1/2} t_1 }{x};q\right)_{n}, \quad c_0=1.
\label{eq:gxg3}
\end{equation}
Set $a_{-1}= a_{-2} =a_{-3} =0$.
The following formulas are shown directly;
\begin{align*}
&\left(x-q^{h_1+1/2}t_1\right)\left(\dfrac{q^{h_1+1/2}t_1}{x/q};q\right)_n = x\left(\dfrac{q^{h_1+1/2}t_1}{x};q\right)_{n+1},  \\
&\left(\dfrac{q^{h_1+1/2}t_1}{qx};q\right)_n =\left(1-\dfrac{q^{h_1+1/2}t_1}{qx}\right)\left(\dfrac{q^{h_1+1/2}t_1}{x};q\right)_{n-1} , \nonumber \\
&\left(1-\dfrac{q^{h_1+1/2}t_1}{x}q^n\right)\left(\dfrac{q^{h_1+1/2}t_1}{x};q\right)_n =\left(\dfrac{q^{h_1+1/2}t_1}{x};q\right)_{n+1} . \nonumber 
\end{align*}
We substitute the function $g(x)$ in Eq.~(\ref{eq:gxg3}) into Eq.~(\ref{eq:qhypervar1}).
Then we have
\begin{align}
&q^{\alpha_1}\displaystyle\sum_{n=0}^{\infty}c_n\left(1-\dfrac{q^{h_2+1/2}t_2}{x}\right)\left(\dfrac{q^{h_1+1/2}t_1}{x};q\right)_{n+1} \label{eq:g3qal1al2}  \\
&-\displaystyle\sum_{n=0}^{\infty}c_n\left\{(q^{\alpha_1}+q^{\alpha_2})+\frac{E}{x}+\dfrac{q^{(h_1+h_2+l_1+l_2+\alpha_1+\alpha_2)/2}(q^{1/2}+q^{-1/2})t_1t_2}{x^2}\right\}\left(\dfrac{q^{h_1+1/2}t_1}{x};q\right)_n \nonumber \\
&+q^{\alpha_2}\displaystyle\sum_{n=0}^{\infty}c_n \left(1-\dfrac{q^{l_1-1/2}t_1}{x}\right)\left(1-\dfrac{q^{l_2-1/2}t_2}{x}\right)\left(1-\dfrac{q^{h_1+1/2}t_1}{qx}\right)\left(\dfrac{q^{h_1+1/2}t_1}{x};q\right)_{n-1}=0 . \nonumber 
\end{align}
Note that $ E= -q^{(h_1 +h_2 + l_1 + l_2 +\alpha _1 +\alpha _2 )/2 } \{ (q^{- h_2 }+q^{-l_2 })t_1 + (q^{- h_1 }+ q^{- l_1 }) t_2 \}  $.
We determine the constants $ A^{[-]}_n$, $B^{[-]}_n$, $A^{[0]}_n$, $B^{[0]}_n$, $C^{[0]}_n$, $A^{[+]}_n$, $B^{[+]}_n$, $C^{[+]}_n$ and $D^{[+]}_n$ by the identities
\begin{align*}
& 1-q^{h_2+1/2}t_2/x=A^{[-]}_n+B^{[-]}_n\left(1-\dfrac{q^{h_1+1/2}t_1}{x}q^{n+1}\right) , \\
& (q^{\alpha_1}+q^{\alpha_2})+\dfrac{E}{x}+\dfrac{q^{(h_1+h_2+l_1+l_2+\alpha_1+\alpha_2)/2}(q^{1/2}+q^{-1/2})t_1t_2}{x^2} \nonumber \\
& = A^{[0]}_n +B^{[0]}_n \left(1-\dfrac{q^{h_1+1/2}t_1}{x}q^n\right)+C^{[0]}_n \left(1-\dfrac{q^{h_1+1/2}t_1}{x}q^n\right)\left(1-\dfrac{q^{h_1+1/2}t_1}{x}q^{n+1}\right) , \nonumber \\
& \left(1-\dfrac{q^{l_1-1/2}t_1}{x}\right)\left(1-\dfrac{q^{l_2-1/2}t_2}{x}\right)\left(1-\dfrac{q^{h_1+1/2}t_1}{qx}\right) \nonumber \\
& =A^{[+]}_n+B^{[+]}_n \left(1-\dfrac{q^{h_1+1/2}t_1}{x}q^{n-1}\right)+C^{[+]}_n \left(1-\dfrac{q^{h_1+1/2}t_1}{x}q^{n-1}\right)\left(1-\dfrac{q^{h_1+1/2}t_1}{x}q^n\right) \nonumber \\
& \qquad \qquad  +D^{[+]}_n \left(1-\dfrac{q^{h_1+1/2}t_1}{x}q^{n-1}\right)\left(1-\dfrac{q^{h_1+1/2}t_1}{x}q^n\right)\left(1-\dfrac{q^{h_1+1/2}t_1}{x}q^{n+1}\right). \nonumber 
\end{align*}
Then Eq.~(\ref{eq:g3qal1al2}) is equivalent to 
\begin{align*}
&0= q^{\alpha_1}\displaystyle \sum_{n=0}^{\infty}c_n\left\{A^{[-]}_n \left(\dfrac{q^{h_1+1/2}t_1}{x};q\right)_{n+1}+B^{[-]}_n \left(\dfrac{q^{h_1+1/2}t_1}{x};q\right)_{n+2}\right\} \\
&-\displaystyle\sum_{n=0}^{\infty}c_n\left\{A^{[0]}_n \left(\dfrac{q^{h_1+1/2}t_1}{x};q\right)_n+B^{[0]}_n \left(\dfrac{q^{h_1+1/2}t_1}{x};q\right)_{n+1}+C^{[0]}_n \left(\dfrac{q^{h_1+1/2}t_1}{x};q\right)_{n+2}\right\}  \nonumber \\
&\; \; +q^{\alpha_2}\displaystyle\sum_{n=0}^{\infty}c_n\left\{A^{[+]}_n \left(\dfrac{q^{h_1+1/2}t_1}{x};q\right)_{n-1}+B^{[+]}_n \left(\dfrac{q^{h_1+1/2}t_1}{x};q\right)_{n} \right.  \nonumber \\
& \qquad \qquad \qquad \left. +C^{[+]}_n \left(\dfrac{q^{h_1+1/2}t_1}{x};q\right)_{n+1}+D^{[+]}_n \left(\dfrac{q^{h_1+1/2}t_1}{x};q\right)_{n+2}\right\}  \nonumber \\
& = \sum_{n=0}^{\infty} Q(c_{n-2},c_{n-1},c_{n},c_{n+1}) \left(\dfrac{q^{h_1+1/2}t_1}{x};q\right)_n ,  \nonumber 
\end{align*}
where 
\begin{align*}
& Q(c_{n-2},c_{n-1},c_{n},c_{n+1})= \{q^{\alpha_1}B^{[-]}_{n-2}+q^{\alpha_2}D^{[+]}_{n-2}-C^{[0]}_{n-2}  \} c_{n-2}\\
& \qquad +\{q^{\alpha_1}A^{[-]}_{n-1} +q^{\alpha_2}C^{[+]}_{n-1}-B^{[0]}_{n-1} \} c_{n-1} +\{q^{\alpha_2}B^{[+]}_{n} -A^{[0]}_{n} \} c_n+q^{\alpha_2} A^{[+]}_{n+1} c_{n+1} . \nonumber
\end{align*}
Set $X=-q^{h_1 -l_2} t_1/t_2 $, $Y_1 = q^{h_1-l_1 } $, $Y_2 = q^{\lambda + \alpha _1 } $ , $Y_3 = q^{h_2-l_2 } $. Then we have
\begin{align*}
& \{ q^{\alpha_1}B^{[-]}_{n-2} +q^{\alpha_2}D^{[+]}_{n-2} - C^{[0]}_{n-2} \} t_1/t_2 \\
&=q^{h_2-h_1-2\lambda -\alpha_1-3n+4}(1-q^{\lambda +\alpha_1+n-1})(1-q^{\lambda +\alpha_1+n-2}) , \nonumber \\
& \{ q^{\alpha_1}A^{[-]}_{n-1}+q^{\alpha_2}C^{[+]}_{n-1} -B^{[0]}_{n-1} \} t_1/t_2 \nonumber \\
&=-q^{h_2-h_1-2\lambda -\alpha_1-2n+2}(1-q^{\lambda +\alpha_1+n-1})\{(1-Y_2Y_3q^{n-1})X \nonumber \\
& \qquad \qquad \qquad -(1-q^{-n+1}-q^{-n}-q^{-n-1}+Y_1+Y_2q^{-1})\} , \nonumber \\
& \{ q^{\alpha_2}B^{[+]}_{n} -A^{[0]}_{n} \} t_1/t_2 \nonumber \\
&=q^{h_2-h_1-2\lambda -\alpha_1-2n}\{(1+q-q^{n+1}-Y_1q^{n+1}+Y_1q^{2n+1}-Y_2q^n-Y_2Y_3q^n+Y_2^2Y_3q^{2n})X \nonumber \\
&-(1+q-q^{-n-1}-q^{-n}-q^{-n+1}+Y_1+Y_1q-Y_1q^{n+1}+Y_2+Y_2q^{-1}-Y_2q^n-Y_1Y_2q^n)\} , \nonumber \\
& q^{\alpha_2}A^{[+]}_{n+1} t_1/t_2 =q^{h_2-h_1-2\lambda -\alpha_1-2n-1}(1-Y_1q^{n+1})(1+Xq^{n+1})(1-q^{-n-1}) . \nonumber 
\end{align*}
Write $c_n=\dfrac{q^n(q^{\lambda +\alpha_1};q)_n}{(q^{h_1-l_2+1}t_1/t_2;q)_n}\tilde{c}_n$.
Then the equation $ Q(c_{n-2},c_{n-1},c_{n},c_{n+1})= 0 $ is equivalent to
\begin{align*}
&q^{-n+2}(1+Xq^{n-1})(1+Xq^n)\tilde{c}_{n-2} \\
&-q\{(1-Y_2Y_3q^{n-1})X-(1-q^{-n+1}-q^{-n}-q^{-n-1}+Y_1+Y_2q^{-1})\}(1+Xq^n)\tilde{c}_{n-1} \nonumber \\
&+\{(1+q-q^{n+1}-Y_1q^{n+1}+Y_1q^{2n+1}-Y_2q^n-Y_2Y_3q^n+Y_2^2Y_3q^{2n})X \nonumber \\
&-(1+q-q^{-n-1}-q^{-n}-q^{-n+1}+Y_1+Y_1q-Y_1q^{n+1}+Y_2+Y_2q^{-1}-Y_2q^n-Y_1Y_2q^n)\}\tilde{c}_n \nonumber \\
&+(1-Y_1q^{n+1})(1-Y_2q^n)(1-q^{-n-1})\tilde{c}_{n+1}=0 . \nonumber 
\end{align*}
By substituting $\displaystyle \tilde{c}_n= \sum_{k=0}^{n}c_{n,k}X^k$, we have
\begin{align*}
&\displaystyle 0= \sum_{k=0}^{n}\{q^{n+1}X^{k+2}+q(1+q)X^{k+1}+q^{-n+2}X^k\}c_{n-2,k} \\
&-\displaystyle\sum_{k=0}^{n}q\{(1-Y_2Y_3q^{n-1})X-(1-q^{-n+1}-q^{-n}-q^{-n-1}+Y_1+Y_2q^{-1})\}(1+Xq^n)X^kc_{n-1,k} \nonumber \\
&+\displaystyle\sum_{k=0}^{n}\{(1+q-q^{n+1}-Y_1q^{n+1}+Y_1q^{2n+1}-Y_2q^n-Y_2Y_3q^n+Y_2^2Y_3q^{2n})X \nonumber \\
&-(1+q-q^{-n-1}-q^{-n}-q^{-n+1}+Y_1+Y_1q-Y_1q^{n+1}+Y_2+Y_2q^{-1}-Y_2q^n-Y_1Y_2q^n)\}X^kc_{n,k} \nonumber \\
&+\displaystyle\sum_{k=0}^{n}(1-Y_1q^{n+1})(1-Y_2q^n)(1-q^{-n-1})X^kc_{n+1,k} \nonumber \\
& =  \sum_{k=0}^{n} \tilde{Q} (n,k ) X^k , \nonumber 
\end{align*}
where 
\begin{align*}
& \tilde{Q} (n,k ) = \\
& q^{n+1} c_{n-2,k-2}+q(1+q)c_{n-2,k-1}+q^{-n+2}c_{n-2,k}-q^{n+1}(1-Y_2Y_3q^{n-1})c_{n-1,k-2} \nonumber \\
& -q(2+q^{-1}+q-q^n-Y_1 q^n-Y_2q^{n-1}-Y_2Y_3q^{n-1})c_{n-1,k-1} \nonumber \\
& +q(1-q^{-n+1}-q^{-n}-q^{-n-1}+Y_1+Y_2q^{-1})c_{n-1,k} \nonumber \\
& +\{ 1+q-q^{n+1}-Y_1(q^{n+1} -q^{2n+1}) -Y_2q^n-Y_2Y_3q^n+Y_2^2Y_3q^{2n} \} c_{n,k-1} \nonumber \\
& -\{ 1+q-q^{-n-1}-q^{-n}-q^{-n+1}+Y_1(1+ q- q^{n+1}) +Y_2(1 + q^{-1}- q^n) -Y_1 Y_2q^n \} c_{n,k} \nonumber \\
& +(1-Y_1q^{n+1})(1-Y_2q^n)(1-q^{-n-1})c_{n+1,k}.  \nonumber 
\end{align*}
We show that 
\begin{align}
& c_{n,k}=q^{k(k+1)/2}\dfrac{(q^{\lambda +\alpha_1+l_2-h_2};q)_k}{(q^{h_1-l_1+1};q)_k(q;q)_k(q;q)_{n-k}}=q^{k(k+1)/2}\dfrac{(Y_2Y_3;q)_k}{(Y_1q;q)_k(q;q)_k(q;q)_{n-k}}
\label{eq:cnk}
\end{align}
satisfies the equation $\tilde{Q} (n,k ) = 0$ for $0\leq k \leq n$.
It follows from $c_{n-2,k}=(1-q^{n-k-1})(1-q^{n-k})c_{n,k}$, $c_{n-1,k}=(1-q^{n-k})c_{n,k}$ and $c_{n+1,k}=\dfrac{1}{1-q^{n-k+1}}c_{n,k}$ that the equation $\tilde{Q} ( n,k ) = 0$ is equivalent to
\begin{align}
&-q^{2k}(1-q^{n-k+1})(1-q^{n-k+2})(1-Y_2Y_3q^{k-2})c_{n,k-2} \label{eq:cnksimp} \\
& +q^{k+1}(1-q^{n-k+1})(1-Y_1q^{k-1})(1-q^{k-1})c_{n,k-1}  \nonumber \\
&+q^k(1-Y_2q^{k-1})(1-Y_2Y_3q^{k-1})(1-q^{n-k+1})c_{n,k-1}  \nonumber \\
& -(1-Y_2q^{k-1})(1-Y_1q^k)(1-q^k)c_{n,k}=0 . \nonumber 
\end{align}
On the other hand, it follows from Eq.~(\ref{eq:cnk}) that 
\begin{align*}
& c_{n,k}=q^k\dfrac{(1-Y_2Y_3q^{k-1})(1-q^{n-k+1})}{(1-Y_1q^k)(1-q^k)}c_{n,k-1} .
\end{align*}
Therefore we obtain Eq.~(\ref{eq:cnksimp}), and the function 
\begin{equation*}
g(x)=x^{-\alpha_1}\sum^{\infty}_{n=0} \left(\frac{q^{h_1+1/2} t_1 }{x};q\right)_{n} \dfrac{q^n(q^{\lambda_1+\alpha_1};q)_n}{(q^{h_1-l_2+1}t_1/t_2;q)_n} \sum_{k=0}^{n}c_{n,k}X^k 
\end{equation*}
is a solution of the variant of the $q$-hypergeometric equation of degree two.
\end{prf}

\section{Several limits of variants of the $q$-hypergeometric equation} \label{sec:limit}
We investigate relationship among the $q$-hypergeometric equation and its variants, and we consider the limit to obtain the differential equations as $q\to 1$.
\subsection{From the variant of the $q$-hypergeometric equation of degree three to the variant of the $q$-hypergeometric equation of degree two}
Recall that the variant of the $q$-hypergeometric equation of degree three was given by
\begin{align*}
& (x-q^{h_1 +1/2} t_1) (x- q^{h_2 +1/2} t_2) (x- q^{h_3 +1/2} t_3) g(x/q) \\
&  + q^{2\alpha  +1} (x - q^{l_1-1/2}t_1 ) (x - q^{l_2 -1/2} t_2) (x - q^{l_3 -1/2} t_3) g(qx) \nonumber \\
& + q^{\alpha } [ - (q + 1 ) x^3 + q^{1/2} \{ (q^{h_1} + q^{l_1 })t_1 + (q^{h_2} + q^{l_2 })t_2 + (q^{h_3} + q^{l_3 })t_3 \} x^2 \nonumber \\
& \quad - q^{(h_1+h_2+h_3+l_1+l_2+l_3 +1)/2} \{ (q^{- h_1 }+q^{-l_1 })t_2 t_3 + (q^{- h_2 }+ q^{- l_2 }) t_1 t_3 \nonumber \\
& \quad + (q^{- h_3 }+ q^{- l_3 }) t_1 t_2  \} x + q^{(h_1 +h_2 + h_3 + l_1 + l_2 + l_3 )/2 } ( q + 1 ) t_1 t_2t_3 ] g(x) =0. \nonumber
\end{align*}
We divide it by $- q^{h_3 +1/2} t_3 $ and consider the limit $t_3 \to \infty $ (i.e.~$1/ t_3 \to 0$) while the variable $x$ is fixed.
Then the equation tends to
\begin{align*}
& (x-q^{h_1 +1/2} t_1) (x- q^{h_2 +1/2} t_2) g(x/q) \\
& + q^{2\alpha + l_3 -h_3 } (x - q^{l_1-1/2}t_1 ) (x - q^{l_2 -1/2} t_2)  g(qx) \nonumber \\
& - [ (q^{\alpha } + q^{\alpha -h_3 + l_3 }) x^2 
- q^{(h_1+h_2+l_1+l_2+2\alpha -h_3+ l_3 )/2} \{ (q^{- h_1 }+q^{-l_1 })t_2 \nonumber \\
& \quad + (q^{- h_2 }+ q^{- l_2 }) t_1 \} x + q^{(h_1 +h_2  + l_1 + l_2 +2 \alpha - h_3 + l_3 )/2 } ( q^{1/2} + q^{-1/2} ) t_1 t_2  ] g(x) =0. \nonumber
\end{align*}
Therefore we obtain the variant of the $q$-hypergeometric equation of degree two given in Eq.~(\ref{eq:qhypervar1}) with the parameter $\{ \alpha _1 ,\alpha _2 \} = \{ \alpha , \alpha -h_3 + l_3  \} $.

We investigate the limit of the solutions given in Conjecture \ref{conj:var3} as $t_3 \to \infty $.
Recall that the parameter $\nu = (h_1 +h_2 +h_3 -l_1-l_2-l_3 +1)/2$ (resp.~the parameter $\lambda =  (h_1 +h_2 -l_1-l_2 -\alpha _1 -\alpha _2 +1)/2  $) was used to express the solutions of the variant of the $q$-hypergeometric equation of degree three (resp. degree two).
For simplicity, we consider the case $\alpha _1 =\alpha $ and $\alpha _2 = \alpha -h_3 + l_3 $.
Then we have $\alpha _1 -\alpha _2 = h_3 - l_3$, $\nu = \lambda + \alpha _1$ and $\nu - h_3 + l_3 = \lambda + \alpha _2$.

We pick up the solution of the variant of the $q$-hypergeometric equation of degree three in Conjecture \ref{conj:var3} (i) for $(i,i',i'') = (1,2,3)$.
Then the solution is written as 
\begin{align}
& h(x)= x^{\nu -\alpha } \sum _{n=0}^{\infty} c_n \Big( \frac{x}{q^{l_1-1/2} t_1} ;q \Big)_n  ,
\end{align}
where
\begin{align}
& c_n =\frac{q^n (q^{\nu };q )_n }{(q^{h_2 - l_1 +1} t_2/t_1 ;q )_n (q^{h_3 - l_1 +1} t_3/t_1 ;q )_n  }  \label{eq:cn1123} \\
& \qquad \cdot \sum _{k_2 =0}^n \sum _{k_3 =0}^{n-k_2} c_{n,k_2,k_3} \Big( \frac{- q^{h_2-l_1}t_2}{t_1} \Big)^{k_2} \Big( \frac{- q^{h_3-l_1}t_3 }{t_1 } \Big)^{k_3} , \nonumber \\
& c_{n,k_2,k_3} =  q^{(k_2+k_3)(k_2 +k_3 +1)/2} \frac{(q^{\nu -h_2+l_2 } ;q)_{k_2} (q^{\nu -h_3+l_3 } ;q)_{k_3} }{(q;q)_{n-k_2-k_3} (q;q)_{k_2}(q;q)_{k_3} (q^{h_1-l_1 +1};q)_{k_2 + k_3} } . \nonumber 
\end{align}
On the limit $t_3 \to \infty $, the terms except for $(k_2 , k_3 ) = (0,n) $ in Eq.~(\ref{eq:cn1123}) vanish, and we have
\begin{align*}
& c_n  \to q^n \frac{(q^{\nu };q )_n }{(q^{h_2 - l_1 +1} t_2/t_1 ;q )_n } \frac{ (q^{\nu -h_3+l_3 } ;q)_{n} }{ (q;q)_{n} (q^{h_1-l_1 +1};q)_{n} } .
\end{align*}

Since $\nu = \lambda + \alpha = \lambda + \alpha _1 $ and $\nu - h_3 + l_3 = \lambda + \alpha _2$, we obtain the function $g(x)$ in Theorem \ref{thm:thm2} for $(i,i')=(1,2)$.
We can also obtain the function $g(x)$ in Theorem \ref{thm:thm2} for $(i,i')=(2,1)$ by considering the limit from Conjecture \ref{conj:var3} (i) in the case $(i,i',i'') = (2,1,3)$.

Next we pick up the solution of the variant of the $q$-hypergeometric equation of degree three in Conjecture \ref{conj:var3} (ii) for $(i,i',i'') = (1,2,3)$.
Then the solution is written as 
\begin{align*}
& h(x)= x^{-\alpha } \sum _{n=0}^{\infty} c_n \Big( \frac{q^{h_1+1/2} t_1 }{x}  ;q \Big)_n   ,
\end{align*}
where
\begin{align}
& c_n = \frac{ q^n (q^{\nu };q )_n }{(q^{h_1 - l_2 +1} t_1/t_2 ;q )_n (q^{h_1 - l_3 +1} t_1/t_3 ;q )_n  } \label{eq:cn2123} \\
& \qquad \cdot \sum _{k_2 =0}^n \sum _{k_3 =0}^{n-k_2} c_{n,k_2,k_3} \Big( \frac{- q^{h_1-l_2}t_1}{t_2 } \Big)^{k_2} \Big( \frac{- q^{h_1-l_3}t_1}{t_3 } \Big)^{k_3} , \nonumber  \\
& c_{n,k_2,k_3} =  q^{(k_2+k_3)(k_2 +k_3 +1)/2} \frac{(q^{\nu -h_2+l_2 } ;q)_{k_2} (q^{\nu -h_3+l_3 } ;q)_{k_3} }{(q;q)_{n-k_2-k_3} (q;q)_{k_2}(q;q)_{k_3} (q^{h_1-l_1 +1};q)_{k_2 + k_3} } . \nonumber 
\end{align}
On the limit $t_3 \to \infty $, the terms except for $k_3=0$ in Eq.~(\ref{eq:cn2123}) vanish, and we have
\begin{align*}
& c_n \to q^n \frac{(q^{\nu };q )_n }{(q^{h_1 - l_2 +1} t_1/t_2 ;q )_n  } \\
& \qquad \qquad \cdot \sum _{k_2 =0}^n  q^{k_2 (k_2 +1)/2} \frac{(q^{\nu -h_2+l_2 } ;q)_{k_2} }{(q;q)_{n-k_2} (q;q)_{k_2} (q^{h_1-l_1 +1};q)_{k_2 } } \Big( \frac{- q^{h_1-l_2}t_1 }{t_2 } \Big) ^{k_2} .
\end{align*}
Since $\nu = \lambda + \alpha _1 $, we obtain the function $g(x)$ in Theorem \ref{thm:thm3} for $(i,i')=(1,2)$.
We can also obtain function the $g(x)$ in Theorem \ref{thm:thm3} for $(i,i')=(2,1)$ by considering the limit from Conjecture \ref{conj:var3} (ii) in the case $(i,i',i'') = (2,1,3)$.

We pick up the solution of the variant of the $q$-hypergeometric equation of degree three in Conjecture \ref{conj:var3} (ii) for $(i,i',i'') = (3,1,2)$.
Then the solution is written as 
\begin{align*}
& h(x)= x^{-\alpha } \sum _{n=0}^{\infty} c_n \Big( \frac{q^{h_3+1/2} t_3 }{x}  ;q \Big)_n ,
\end{align*}
whose coefficients are determined by
\begin{align}
& c_n =q^n \frac{(q^{\nu };q )_n }{(q^{h_3 - l_1 +1} t_3/t_1 ;q )_n (q^{h_3 - l_2 +1} t_3/t_2 ;q )_n  } \label{eq:cn2312} \\
& \qquad \cdot \sum _{k_2 =0}^n \sum _{k_1 =0}^{n-k_2} c_{n,k_1,k_2} (- q^{h_3 -l_1}t_3 /t_1)^{k_1} (- q^{h_3 -l_2}t_3 /t_2)^{k_2} , \nonumber \\
& c_{n,k_1,k_2} =  q^{(k_1+k_2)(k_1 +k_2 +1)/2} \frac{(q^{\nu -h_1+l_1 } ;q)_{k_1} (q^{\nu -h_2+l_2 } ;q)_{k_2}  }{(q;q)_{n-k_1-k_2} (q;q)_{k_1}(q;q)_{k_2} (q^{h_3 -l_3 +1};q)_{k_1 + k_2} }. \nonumber 
\end{align}
On the limit $t_3 \to \infty $, the terms except for $k_2+k_3=n$  in Eq.~(\ref{eq:cn2312}) vanish, and we have
\begin{align*}
& c_n \Big( \frac{q^{h_3+1/2} t_3 }{x}  ;q \Big) _n  \\
& \to \Big( \frac{q^{ 1/2}}{x} \Big) ^n \frac{(q^{\nu };q )_n }{ (q^{h_3 -l_3 +1};q)_{n} } \sum _{k =0}^n \frac{(q^{\nu -h_1+l_1 } ;q)_{n-k} (q^{\nu -h_2+l_2 } ;q)_{k} }{ (q;q)_{k}(q;q)_{n-k} } (q^{l_1 } t_1)^{k} ( q^{l_2 } t_2)^{n- k } . \nonumber 
\end{align*}
Since $\alpha =\alpha _1  $ and $\nu = \lambda + \alpha _1 $, we obtain the function $g(x)$ in Theorem \ref{thm:thm1}.

In summary, we recovered the solutions of the variant of the $q$-hypergeometric equation of degree two in Theorems \ref{thm:thm1}, \ref{thm:thm2} and \ref{thm:thm3} by the limit from the solutions of the variant of the $q$-hypergeometric equation of degree three in Conjecture \ref{conj:var3}.

\subsection{From the variant of the $q$-hypergeometric equation of degree two to the $q$-hypergeometric equation}
We investigate the limit from the variant of the $q$-hypergeometric equation of degree two to the $q$-hypergeometric equation.
Recall that the variant of the $q$-hypergeometric equation of degree two was given by
\begin{align*}
& 0= (x-q^{h_1 +1/2} t_1) (x - q^{h_2 +1/2} t_2) g(x/q) \\
& \quad  + q^{\alpha _1 +\alpha _2} (x - q^{l_1-1/2}t_1 ) (x - q^{l_2 -1/2} t_2) g(q x) -[ (q^{\alpha _1} +q^{\alpha _2} ) x^2 \nonumber \\
&  \quad  -p \{ (q^{- h_2 }+q^{-l_2 })t_1 + (q^{- h_1 }+ q^{- l_1 }) t_2 \} x + p ( q^{1/2}+ q^{-1/2}) t_1 t_2 ] g(x) , \nonumber
\end{align*}
where $p= q^{(h_1 +h_2 + l_1 + l_2 +\alpha _1 +\alpha _2 )/2 } $.
We take the limit $t_2 \to 0$.
Then we have
\begin{align}
& (x-q^{h_1 +1/2} t_1) g(x/q)  + q^{\alpha _1 +\alpha _2} (x - q^{l_1-1/2}t_1 ) g(q x) \label{eq:qhypeqnt1} \\
&  \qquad -[ (q^{\alpha _1} +q^{\alpha _2} ) x -q^{(h_1 +h_2 + l_1 + l_2 +\alpha _1 +\alpha _2 )/2 } (q^{- h_2 }+q^{-l_2 })t_1  ] g(x) =0. \nonumber
\end{align}
The standard form of the $q$-hypergeometric equation given in Eq.~(\ref{eq:qhyp}) is realized by setting 
\begin{align}
& t_1=1, \; h_1 =1/2, \; h_2-l_2 = \alpha _1 +\alpha _2 + l_1 - 3/2 , \label{eq:restr} \\
& a= q^{\alpha _1 }  , \; b= q^{\alpha _2 }  , \; c= q^{\alpha _1 +\alpha _2 + l_1-1/2} . \nonumber
\end{align}

We investigate the limit of the solutions given in Theorems \ref{thm:thm1}, \ref{thm:thm2} and \ref{thm:thm3}  as $t_2 \to 0$.
We pick up the solution of the variant of the $q$-hypergeometric equation of degree two in Theorem \ref{thm:thm1} for $(i,i') = (1,2)$.
Then it is written as 
\begin{align*}
g (x) = & x^{-\alpha _1 } \sum _{n=0}^{\infty} (q^{1/2} x^{-1})^n \frac{ ( q^{\lambda +\alpha _1 } ; q )_n }{(q^{\alpha _1  -\alpha _2 +1 } ; q )_n } \\
& \cdot \sum _{k=0}^n  \frac{(q^{ \lambda +\alpha _1  -h_2 +l_2 }; q)_k (q^{  \lambda +\alpha _1  -h_1 +l_1 } ;q)_{n-k }}{(q;q )_k (q;q)_{n-k}}  (q^{l_1} t_1)^k (q ^{l_2}t_2 ) ^{n-k} , \nonumber 
\end{align*}
where $\lambda = (h_1 +h_2 -l_1-l_2 -\alpha _1-\alpha _2+1 )/2 $.
We take the limit $t_2 \to 0$.
Then the terms except for $k=n$ vanish and we have
\begin{align}
g (x) \to 
x^{-\alpha _1 } \sum _{n=0}^{\infty} \frac{ ( q^{\lambda +\alpha _1 } ; q )_n }{(q^{\alpha _1  -\alpha _2 +1 } ; q )_n } \frac{(q^{ \lambda +\alpha _1  -h_2 +l_2 }; q)_n }{(q;q )_n }  \Big( \frac{q^{l_1 +1/2 } t_1 }{x} \Big)^n . \label{eq:g1qhyp}
\end{align}
Then we can show directly that the right-hand side of Eq.~(\ref{eq:g1qhyp}) satisfies Eq.~(\ref{eq:qhypeqnt1}).
Under the relation in Eq.~(\ref{eq:restr}), we have the relation $\lambda =0$ and we recover Eq.~(\ref{eq:qHGSinfty}).

The solution of the variant of the $q$-hypergeometric equation of degree two in Theorem \ref{thm:thm2} for $(i,i') = (1,2)$ is written as
\begin{align*}
& g(x)= x^{\lambda } \sum _{n=0}^{\infty} \Big( \frac{x}{q^{l_1-1/2} t_1 } ;q \Big)_n q^n \frac{(q^{\lambda +\alpha _1 };q )_n (q^{\lambda +\alpha _2 };q )_n }{(q^{h_1 - l_1 +1};q )_n (q^{h_2 - l_1 +1} t_2/t_1 ;q )_n (q;q)_n } . 
\end{align*}
As $t_2 \to 0$, we have
\begin{align*}
& g (x) \to x^{\lambda } \sum _{n=0}^{\infty} \Big( \frac{x}{q^{l_1-1/2} t_1 } ;q \Big)_n  \frac{(q^{\lambda +\alpha _1 };q )_n (q^{\lambda +\alpha _2 };q )_n }{(q^{h_1 - l_1 +1};q )_n  (q;q)_n } q^n,
\end{align*}
and the function obtained by the limit is a solution to Eq.~(\ref{eq:qhypeqnt1}).
Under the relation in Eq.~(\ref{eq:restr}), we recover Eq.~(\ref{eq:qHGSPoch}).

The solution of the variant of the $q$-hypergeometric equation of degree two in Theorem \ref{thm:thm2} for $(i,i') = (2,1)$ is written as
\begin{align*}
& g (x)=  x^{\lambda } \sum _{n=0}^{\infty} \Big( \frac{x}{q^{l_2-1/2} t_2 } ;q \Big)_n q^n \frac{(q^{\lambda +\alpha _1 };q )_n (q^{\lambda +\alpha _2 };q )_n }{(q^{h_2 - l_2 +1};q )_n (q^{h_1 - l_2 +1} t_1/t_2 ;q )_n (q;q)_n } .
\end{align*}
As $t_2 \to 0$, we have
\begin{align*}
 g (x) \to x^{\lambda } \sum _{n=0}^{\infty} \frac{(q^{\lambda +\alpha _1 };q )_n (q^{\lambda +\alpha _2 };q )_n }{(q^{h_2 - l_2 +1};q )_n (q;q)_n } \Big( \frac{x }{q^{h_1 -1/2 } t_1 } \Big)^n ,
\end{align*}
and we obtain a solution to Eq.~(\ref{eq:qhypeqnt1}).
Under the relation in Eq.~(\ref{eq:restr}), we recover Eq.~(\ref{eq:qhypser}).

The solution of the variant of the $q$-hypergeometric equation of degree two in Theorem \ref{thm:thm3} for $(i,i') = (1,2)$ is written as
\begin{align*}
& g (x)= x^{-\alpha _1 } \sum _{n=0}^{\infty} \Big( \frac{q^{h_1+1/2} t_1 }{x} ;q \Big)_n \frac{q^n (q^{\lambda +\alpha _1 };q )_n }{(q^{h_1 - l_2 +1} t_1/t_2 ;q )_n  } \cdot \\
& \qquad \qquad \qquad \qquad \cdot \sum _{k=0}^n   \frac{(q^{\lambda -h_2+l_2+\alpha _1} ;q)_k q^{k(k+1)/2} (-q^{h_1-l_2}t_1/t_2)^k }{(q^{h_1-l_1+1};q)_{k} (q;q)_k (q;q)_{n-k} }  .
\nonumber
\end{align*}
As $t_2 \to 0$, we have
\begin{align*}
 g (x) \to x^{-\alpha _1 } \sum _{n=0}^{\infty} \Big( \frac{q^{h_1+1/2} t_1 }{x}  ;q \Big)_n  \frac{ (q^{\lambda +\alpha _1 };q )_n (q^{\lambda +\alpha _1 -h_2+l_2 } ;q)_n }{(q^{h_1-l_1+1};q)_{n} (q;q)_n  } q^n ,
\end{align*}
and we obtain a solution to Eq.~(\ref{eq:qhypeqnt1}).

\subsection{Limit to the differential equation}
We are going to obtain a differential equation from the difference equation by taking a suitable limit.

Recall that the variant of the $q$-hypergeometric equation of degree two was given by
\begin{align*}
& (x-q^{h_1 +1/2} t_1) (x - q^{h_2 +1/2} t_2) g(x/q) \\
&  \quad  + q^{\alpha _1 +\alpha _2} (x - q^{l_1-1/2}t_1 ) (x - q^{l_2 -1/2} t_2) g(q x) \nonumber \\
&  \quad -[ (q^{\alpha _1} +q^{\alpha _2} ) x^2 +E x + p ( q^{1/2}+ q^{-1/2}) t_1 t_2 ] g(x) =0, \nonumber \\
& p= q^{(h_1 +h_2 + l_1 + l_2 +\alpha _1 +\alpha _2 )/2 } , \quad E= -p \{ (q^{- h_2 }+q^{-l_2 })t_1 + (q^{- h_1 }+ q^{- l_1 }) t_2 \} . \nonumber 
\end{align*}
Set $q=1+ \varepsilon $.
By using Taylor's expansion
\begin{align*}
& g(x/q) =g(x) + (-\varepsilon +\varepsilon ^2 )xg'(x) + \varepsilon ^2 x^2 g''(x) /2 +O( \varepsilon ^3),\\
& g(qx) =g(x) + \varepsilon  x g'(x) + \varepsilon ^2 x^2 g''(x) /2 +O( \varepsilon ^3), \nonumber 
\end{align*}
we find the following limit as $\varepsilon  \to 0$;
\begin{align*}
& x^2(x-t_1)(x-t_2) g''(x) \\
& +  [ (1+ h_2 - l_2)x(x- t_1 ) + (1 +h_1 - l_1 )x(x- t_2) -2\lambda (x -t_1)(x-t_2) ]  xg'(x) \nonumber \\
& + [ \alpha _1 \alpha _2 x^2 + \tilde{B} x + t_1 t_2 \lambda (\lambda +1)] g(x) =0, \nonumber 
\end{align*}
where
\begin{align*}
& \lambda = (h_1+h_2-l_1-l_2 -\alpha _1 -\alpha _2 +1)/2 , \\
& \tilde{B} = -\lambda (\lambda -h_2+l_2 )  t_1 -\lambda (\lambda -h_1+l_1 )  t_2  .  \nonumber 
\end{align*}
This is a Fuchsian differential equation with four singularities $\{ 0,t_1, t_2, \infty \}$ and the local exponents are given by the following Riemann scheme;
\begin{equation*}
\begin{pmatrix}
x=0 &  x=t_1 &  x=t_2 &  x=\infty  \\
\lambda & 0 & 0 & \alpha _1 \\ 
\lambda +1 & l_1-h_1 & l_2-h_2 & \alpha _2 
\end{pmatrix} .
\end{equation*}
Set $y =x^{-\lambda } g(x)$.
Then the singularity $x=0 $ disappears and we have
\begin{align*}
& (x-t_1) (x-t_2) \frac{d^2y}{dx^2} + \{ (2+h_1+h_2-l_1-l_2)(x-t_1 ) \\
& \quad + (t_1-t_2)(h_1 -l_1+1) \}  \frac{dy}{dx} + ( \lambda +\alpha _1) (\lambda +\alpha _2) y= 0. \nonumber 
\end{align*}
The hypergeometric equation appears by setting $z= (x-t_1)/(t_2-t_1)  $,~i.~e.
\begin{align*}
& z(1-z) \frac{d^2y}{dz^2} - \{ (2+h_1+h_2-l_1-l_2)z+l_1-h_1-1 \}  \frac{dy}{dz} \\
&  \quad - ( \lambda +\alpha _1) (\lambda +\alpha _2) y= 0. \nonumber 
\end{align*}

Recall that the variant of the $q$-hypergeometric equation of degree three was given by
\begin{align*}
& (x-q^{h_1 +1/2} t_1) (x- q^{h_2 +1/2} t_2) (x- q^{h_3 +1/2} t_3) g(x/q) \\
&  + q^{2\alpha  +1} (x - q^{l_1-1/2}t_1 ) (x - q^{l_2 -1/2} t_2) (x - q^{l_3 -1/2} t_3) g(qx) \nonumber \\
& + q^{\alpha } [ - (q + 1 ) x^3 + q^{1/2} \{ (q^{h_1} + q^{l_1 })t_1 + (q^{h_2} + q^{l_2 })t_2 + (q^{h_3} + q^{l_3 })t_3 \} x^2 \nonumber \\
& \quad - q^{(h_1+h_2+h_3+l_1+l_2+l_3 +1)/2} \{ (q^{- h_1 }+q^{-l_1 })t_2 t_3 + (q^{- h_2 }+ q^{- l_2 }) t_1 t_3 \nonumber \\
& \quad  + (q^{- h_3 }+ q^{- l_3 }) t_1 t_2  \} x + q^{(h_1 +h_2 + h_3 + l_1 + l_2 + l_3 )/2 } ( q + 1 ) t_1 t_2t_3 ] g(x) =0. \nonumber
\end{align*}
Set $q=1+ \varepsilon $.
By using Taylor's expansion, we find the following limit as $\varepsilon  \to 0$;
\begin{align*}
& x^2(x-t_1)(x-t_2) (x-t_3 )  g''(x)  \\
& + x^2(x-t_1)(x-t_2) (x-t_3 ) \Big\{ \frac{h_1 - l_1 +1 }{x- t_1} + \frac{h_2 - l_2 +1 }{x- t_2} \nonumber \\
& \qquad \qquad +\frac{h_3 - l_3 +1 }{x- t_3} - \frac{ 2( \nu - \alpha ) }{x} \Big\} g'(x) \nonumber \\
& +  [ \alpha (\alpha +1 ) x^3 - \alpha \{ (l_1 -h_1 +\alpha )t_1+(l_2 -h_2 +\alpha )t_2 \nonumber \\
& \qquad \qquad   +(l_3 -h_3 +\alpha ) t_3 \} x^2 + \tilde{B} x - t_1 t_2 t_3( \nu  -\alpha  ) (\nu -\alpha + 1) ]g(x) =0,  \nonumber 
\end{align*}
where $\nu = ( h_1 +h_2 +h_3 -l_1 -l_2 -l_3 +1 )/2$ and  
\begin{align*}
& \tilde{B} = (\nu -\alpha +1/2 )\{ (\nu -\alpha - h_1 + l_1 ) t_2 t_3 +( \nu -\alpha - h_2 + l_2  ) t_3 t_1 \\
& \qquad \qquad \qquad + ( \nu -\alpha - h_3 + l_3  )t_1 t_2 \}. \nonumber
\end{align*}
By setting $y =x^{-\nu +\alpha }(x-t_2)^{\nu } g(x)$ and $z= (x-t_1)(t_3-t_2)/\{(x-t_2)(t_3-t_1) \} $, we obtain the hypergeometric equation
\begin{align*}
& z(1-z) \frac{d^2y}{dz^2} - \{ (2+h_1+h_3-l_1-l_3)z+l_1-h_1-1 \}  \frac{dy}{dz} \\
& \qquad \qquad - \nu ( \nu -h_2 +l_2 ) y= 0. \nonumber 
\end{align*}

\section{Concluding remarks} \label{sec:concl}

In this paper, we introduced the variants of the $q$-hypergeometric equation of degree two or three, which corresponds to the Fuchsian differential equation with three singularities which may not contain $0$ or/and $\infty $ as $q\to 1$.
We obtained three kinds of explicit solutions to the variant of the $q$-hypergeometric equation of degree two, and also obtained the conjectural solutions to the variant of the $q$-hypergeometric equation of degree three.

It is known that the relations between the solutions of the hypergeometric differential equation about the origin and those about the point at infinity are written explicitly by using the gamma function, and also the analogous relations are known for the basic hypergeometric functions (see \cite{GR}).
The integral representations of the hypergeometric functions or the basic hypergeometric functions were applied to obtain these relations.

Then we propose a problem to obtain the relations of three kinds of solutions to the variants of the $q$-hypergeometric equation of degree two.
For this purpose, integral representations of solutions to the variants of the $q$-hypergeometric equation should be studied well.
Note that the Jackson integral representation of the $q$-Appell function was written in the book of Gasper and Rahman \cite{GR}, and it is possible to obtain a Jackson integral representation of the function $g_1 (x)$ in the introduction by restricting the parameters as in the proof of Theorem \ref{thm:thm1}.

Yamaguchi \cite{Yg} studied the Jackson integral representation of solutions to $q$-differential equations which are related to the variants of the $q$-hypergeometric equation in his master thesis.
Sakai and Yamaguchi \cite{SY} constructed a $q$-analogue of the middle convolution, and it would be useful to study integral representation of solutions to the variants of the $q$-hypergeometric equation or the $q$-Heun equation.

One of our motivations for introducing the variants of the $q$-hypergeometric equation is to understand the confluence of singularities of the $q$-difference equations.
We hope to investigate the procedures of the confluence from the variants of the $q$-hypergeometric equation and their explicit solutions, which has been partly done in \cite{MST}.

\section*{Acknowledgements}
The fourth author was supported by JSPS KAKENHI Grant Number JP18K03378.

\appendix
\section{The $q$-Appell function and a variant of the $q$-hypergeometric equation} \label{sec:qAppell}
It is known in \cite{GR} that the $q$-Appell series defined by 
\begin{equation}
f(x,y)= \Phi ^{(1)}(a;b,b';c;q;x,y) =\sum_{m=0}^\infty \sum_{n=0}^\infty \frac{(a;q)_{m+n}(b;q)_m (b';q)_n}{(c;q)_{m+n}(q;q)_m(q;q)_n}x^my^n
\label{q-App}
\end{equation}
satisfies the difference equations
\begin{align}
& (abx -c/q) f(q^2x,qy) - (bx -1)f(qx,y) \label{a} \\
& \qquad -(ax - c/q) f(qx,qy) +(x-1)f(x,y) = 0, \nonumber \\ 
& (ab'y -c/q )f(qx,q^2y) - (b'y -1)f(x,qy) \nonumber \\
& \qquad -(ay - c/q ) f(qx,qy) +(y-1)f(x,y) =0, \nonumber
\end{align}
which are confirmed by a direct calculation.
By replacing $x$ with $qx$ in the second equation of Eq.~(\ref{a}), we have 
\begin{align}
& (ab'y - c/q)f(q^2x,q^2y) -(b'y -1)f(qx,qy) \label{d} \\
& \qquad - (ay -c/q )f(q^2x,qy)+(y-1)f(qx,y) =0 .  \nonumber
\end{align}
By erasing the term $f(q^2x,qy) $ in Eqs.~(\ref{a}) and (\ref{d}), we have
\begin{align}
& (a- c/q)(bx-y)f(qx,y) = (ay -c/q)(x-1)f(x,y) \label{e} \\
& \qquad -\left\{ (ay -c/q)(ax -c/q) + (abx- c/q)(b'y -1)\right\}f(qx,qy) \nonumber \\
& \qquad +(abx- c/q)(ab'y- c/q)f(q^2x,q^2y) . \nonumber
\end{align}
By erasing the term $f(qx,y) $ in Eqs.~(\ref{a}) and (\ref{d}), we have
\begin{align}
& -(a-c/q)(bx-y)f(q^2x,qy) = -(bx-1)(ab'y-c/q)f(q^2x,q^2y) \label{f} \\
& +\left\{ (y-1)(ax-c/q) + (b'y-1)(bx-1 )\right\} f(qx,qy)  -(x-1)(y-1)f(x,y) .\nonumber 
\end{align}
We replace $(x,y) $ with $(qx,qy)$ in Eq.~(\ref{e}) and substitute it into Eq.~(\ref{f}).
Then we obtain the following third order $q$-difference equation.
\begin{align*}
&  (abqx-c/q)(ab' qy-c/q)f(q^3 x,q^3y) \\
& - \{ (aqx -c/q)(aqy -c/q ) + (abqx-c/q)(b'qy -1)  \nonumber  \\
& \qquad + (bx -1)(ab' q y-c)\} f(q^2x,q^2y) \nonumber  \\
& + \{ (qx-1)(aqy -c/q) + (aqx -c)(y-1) + q(bx-1)(b'y -1 ) \} f(qx,qy) \nonumber \\
& - q(x-1)(y-1)f(x,y) =0 .\nonumber 
\end{align*}
In the case $c= bb'$, the third order $q$-difference equation is reducible and we obtain the following second order $q$-difference equation;
\begin{prop} \label{prop:qApp2ndorder}
The $q$-Appell series $f(x,y)$ defined in Eq.~(\ref{q-App}) with the condition $c= bb'$ satisfies 
\begin{align}
&  (aqx-b')(aqy-b)f(q^2 x,q^2y) \label{2q-Appell} \\
&  -\{ aq(q+1) xy - q(a+b)x-q(a+b')y +(bb' +q) \} f(qx,qy) \nonumber \\
&  +q(x-1)(y-1)f(x,y) =0 . \nonumber
\end{align}
\end{prop}
\begin{prf}
By substituting $f(x,y)=x^m y^n $ into the left-hand side of Eq.~(\ref{2q-Appell}), we have 
\begin{align*}
& q(1-aq^{m+n+1})(1-aq^{m+n}) x^{m+1} y^{n+1} - q(1-aq^{m+n})(1-bq^{m+n}) x^{m+1} y^{n} \\
&  -q(1-aq^{m+n})(1-b'q^{m+n}) x^m y^{n+1} +q(1-cq^{m+n-1})(1-q^{m+n}) x^m y^n . \nonumber
\end{align*}
We substitute the function $f(x,y)$ in Eq.~(\ref{q-App}) with the condition $c = bb'$ into the left-hand side of Eq.~(\ref{2q-Appell}).
It follows from the formula $(1-aq^{m+n+1})(1-aq^{m+n}) (a;q )_{m+n} = (a;q )_{m+n+2} $ and so on that the left-hand side of Eq.~(\ref{2q-Appell}) is equal to  
\begin{align}
&q\left\{ \sum _{m\geq 0, n \geq 0} \frac{(a;q)_{m+n+2}(b;q)_m(b';q)_n}{(c;q)_{m+n}(q;q)_m(q;q)_n}x^{m+1}y^{n+1} \right. \label{3q-Appell} \\
& -\sum _{m\geq 0, n \geq 0}\frac{(a;q)_{m+n+1}(b;q)_m(b';q)_n}{(c;q)_{m+n}(q;q)_m(q;q)_n}(1-bq^{m+n})x^{m+1}y^n \nonumber \\
& -\sum _{m\geq 0, n \geq 0}\frac{(a;q)_{m+n+1}(b;q)_m(b';q)_n}{(c;q)_{m+n}(q;q)_m(q;q)_n}(1-b'q^{m+n})x^my^{n+1} \nonumber \\
& \left. +\sum _{m\geq 0, n \geq 0} \frac{(a;q)_{m+n}(b;q)_m(b';q)_n}{(c;q)_{m+n-1}(q;q)_m(q;q)_n}(1-q^{m+n})x^my^n \right\} \nonumber \\
 =& q \left\{ \sum _{m\geq 1, n \geq 1} \frac{(a;q)_{m+n}(b;q)_{m-1} (b';q)_{n-1}}{(c;q)_{m+n-2}(q;q)_{m-1}(q;q)_{n-1}}x^my^n \right. \nonumber \\
& -\sum _{m\geq 1, n \geq 0} \frac{(a;q)_{m+n}(b;q)_{m-1}(b';q)_{n}}{(c;q)_{m+n-1}(q;q)_{m-1}(q;q)_{n}}(1-bq^{m+n-1})x^my^n \nonumber \\
& -\sum _{m\geq 0, n \geq 1} \frac{(a;q)_{m+n}(b;q)_m(b';q)_{n-1}}{(c;q)_{m+n-1}(q;q)_m(q;q)_{n-1}}(1-b'q^{m+n-1})x^my^n \nonumber \\
& \left. + \sum _{m\geq 0, n \geq 0} \frac{(a;q)_{m+n}(b;q)_m(b';q)_n}{(c;q)_{m+n-1}(q;q)_m(q;q)_n}(1-q^{m+n})x^my^n \right\} . \nonumber 
\end{align}
In the summation, we divide into the case $m\geq 1, n \geq 1 $ and the other cases.
Then the right-hand side of Eq.~(\ref{3q-Appell}) is equal to 
\begin{align*}
& q \left[ \sum^{\infty}_{m=1} \sum^{\infty}_{n=1} \frac{(a;q)_{m+n}(b;q)_{m-1}(b';q)_{n-1}}{(c;q)_{m+n-1}(q;q)_{m}(q;q)_{n}}
\{(1-cq^{m+n-2})(1-q^m)(1-q^n) \right. \\ 
& \quad -(1-b'q^{n-1})(1-q^m)(1-bq^{m+n-1})  -(1-bq^{m-1})(1-q^n)(1-b'q^{m+n-1}) \nonumber \\
& \quad +(1-bq^{m-1})(1-b'q^{n-1})(1-q^{m+n}) \}x^my^n \nonumber \\
& - \sum_{m=1}^\infty \frac{(a;q)_m (b;q)_{m-1}}{(c;q)_{m-1}(q;q)_{m-1}}(1-bq^{m-1})x^m 
-\sum _{n=1}^\infty \frac{(a;q)_n(b';q)_{n-1}}{(c;q)_{n-1}(q;q)_{n-1}}(1-b'q^{n-1})y^n \nonumber \\
& \left. +\sum_{m=1}^\infty \frac{(a;q)_m(b;q)_m}{(c;q)_{m-1}(q;q)_m}(1-q^m)x^m 
+\sum_{n=1}^\infty \frac{(a;q)_n(b';q)_n}{(c;q)_{n-1}(q;q)_n}(1-q^n)y^n \right] , \nonumber 
\end{align*}
and it is equal to $0$, because the second summation is equal to the fourth summation, the third summation is equal to the fifth summation, and  the inside of the curly brace in the first summand is equal to $0$ by the condition $c=b b' $.
Therefore we obtain Eq.~(\ref{2q-Appell}) for the $q$-Appell series with the condition $c= bb'$.
\end{prf}

\end{document}